\documentclass[fleqn,10pt]{wlscirep}
\usepackage[utf8]{inputenc}
\usepackage[T1]{fontenc}

\usepackage{amssymb,amsmath,color,amsfonts}
\usepackage{graphicx}
\graphicspath{{./IMG/}}
\usepackage{caption}
\usepackage{graphicx} 
\usepackage{array}
\usepackage{mathtools}

\usepackage{subfigure}
\usepackage{hyperref}
\usepackage{url}

\newcommand{\mc}{\mathcal}

\newcommand{\real}{\mathbb{R}}

\newcommand{\realpos}{\mathbb{R}_{\geq 0}}

\newcommand\oprocendsymbol{\hbox{$\square$}}
\newcommand\oprocend{\relax\ifmmode\else\unskip\hfill%
\fi\oprocendsymbol}
\renewcommand{\figurename}{Fig.}

\newcommand{\map}[3]{#1: #2 \rightarrow #3}
\newcommand{\setdef}[2]{\{#1 \; : \; #2\}}

\newcommand{\sbs}[2]{{#1}_{\textup{#2}}}

\newcommand{\until}[1]{\{1,\dots,#1\}}

\newcommand{\norm}[1]{\Vert #1 \Vert}

\renewcommand{\i}{\imath}

\usepackage{epstopdf}
\title{Planning a Return to Normal after the COVID-19 Pandemic: Identifying Safe Contact 
Levels via Online Optimization}

\author[1,*]{Gianluca Bianchin}
\author[1]{Emiliano Dall'Anese}
\author[1]{Jorge I. Poveda}
\author[2]{David Jacobson}
\author[3]{Elizabeth J. Carlton}
\author[4]{Andrea G. Buchwald}
\affil[1]{Department of Electrical, Computer, and Energy Engineering, University of Colorado Boulder, CO, USA}
\affil[2]{VanaData, Boulder, CO, USA}
\affil[3]{Department of Environmental and Occupational Health, Colorado School of Public Health, University of Colorado Anschutz, CO, USA}
\affil[4]{Department of Biostatistics and Informatics, 
Colorado School of Public Health, University of Colorado Anschutz, CO, USA}

\affil[*]{Corresponding author -- gianluca.bianchin@colorado.edu}

\begin{abstract}
Since the early months of 2020, non-pharmaceutical interventions (NPIs) --  implemented at varying levels of severity and based on widely-divergent perspectives of risk tolerance -- have been the primary means to control SARS-CoV-2 transmission. We seek to identify how risk tolerance and vaccination rates impact
the rate at which a population can return to pre-pandemic contact behavior. To this end, we develop a novel feedback control method for data-driven decision-making to identify optimal levels of NPIs across geographical regions in order to guarantee that hospitalizations will not exceed a given risk tolerance. Results are shown for the state of Colorado, and they suggest that: coordination in decision-making across regions is essential to maintain the daily number of hospitalizations below the desired limits; increasing risk tolerance can decrease the number of days required to discontinue NPIs, at the cost of an increased number of deaths; and if vaccination uptake is less than $70\%$, at most levels of risk tolerance, return to pre-pandemic contact behaviors
before the early months of 2022 may newly jeopardize the healthcare system.
\end{abstract}
\begin{document}

\flushbottom
\maketitle
%
%
\thispagestyle{empty}

\section*{Introduction}

The primary strategy for mitigating the spread of SARS-CoV-2, to date, has relied on the use of non-pharmaceutical 
interventions (NPIs). Common NPIs include, at various levels of severity, lockdowns, travel restrictions, contact tracing, mask 
wearing, and individual behavioral change.  
These policy-based restrictions and individual behavior changes have had wide-ranging social consequences, including 
disruptive impacts on economies, and they have severely affected the well-being of families and children due to confinement stress and social
disruptions \cite{PO-TA:21,HP-MW-DB:20}. On December 14, 2020, a mass vaccination campaign was initiated in the 
United States and, as the spread of the SARS-CoV-2 virus slowed, 
individuals and policy-makers alike are planning 
a ``return to normality'', where policy makers begin to lift NPIs and individual social 
behaviors are tentatively resumed. The individual and policy-level decisions to lift NPIs and ``return to normal'' depend on: (i) the trajectory of the SARS-CoV-2 states in the absence of NPIs (are immunity levels from vaccination and prior infection enough to prevent future waves of infection?) and (ii) the risk tolerance of individuals and policy-makers. Whether or not this return to normal is safe  is a question of (i) population level immunity (are immunity levels -- either due to previous infection or vaccination -- sufficient to prevent future waves of infection?) and (ii) 
dependent on the risk tolerance of both individuals and 
policy-makers. 

Predicting the minimum level of immunity that is needed to remove all NPIs and return to normal behavior remains an open 
research question, especially in the presence of regional discrepancies in risk tolerance and vaccination rates. 
Tolerable levels of transmission risk have been a controversial topic throughout the pandemic, and could be defined based on numbers of cases, hospitalizations, or deaths. We chose hospitalizations for this analysis because they are a marker of severe COVID-19 disease and less lagged than deaths. In these scenarios, we are not pursuing a goal of elimination – which, with evidence of waning immunity and continued global import, is improbable in the near-term. We instead focus on staying below a tolerable level of endemic disease burden. Defining the tolerable risk level is outside the scope of this paper – but may be based on infrastructure limits such as hospital capacity, as well as societal values including equity, a desire to minimize risk of death or severe disease, or a desire to minimize disruption (e.g., school closures). Here we examine a range of possible risk tolerance thresholds as a demonstration.

In this work, we aim to answer the following questions: given the ongoing vaccination campaigns and a range of levels of risk tolerance, when can 
most NPIs be relaxed (we use the term NPI to encompass all policy restrictions as well as individual-level behavior change) and life return to normal? How is this return to normal impacted by variations in  vaccination rates and risk tolerance? And, how does mobility between regions (counties, states, etc)  impact the return to normal and how to coordinate NPIs across regions?

To answer these questions, we adopt a novel approach for data-driven policy-making using multidisciplinary methods pulling from transmission modeling~\cite{WOK-AGM:27,EF:20,IZK-JCM-PLS:17, moore2021vaccination}, online optimization of dynamical systems~\cite{brunner2012feedback,MC-ED-AB:19,GB-JC-JP-ED:21-tcns,hauswirth2020timescale}, and feedback  control~\cite{khalil2002nonlinear}. We develop a method that allows us to identify maximum contact levels across geographical regions that guarantee hospitalizations will not exceed a given threshold, while minimizing the economic and social  impacts (see Discussion and Methods). This aspect is of 
utmost importance to ensure that future waves of infections do not threaten the stability of public health  infrastructures. The operating principles of the proposed method depart from standard approaches in epidemic control that may be based on specific models for the evolution of the epidemic, and may rely on ad hoc decision rules; the proposed controller operates in closed-loop, in the sense that it suggests whether the level of NPIs should increase (or can be decreased) to meet predefined economic objectives and risk tolerance metrics based on the current level of infections and a prediction of the peak of hospitalizations.

The approach is widely applicable and can be implemented at various geographical granularity and to other disease systems. As a test case, we considered the state of Colorado, USA. We calibrated our models using real-world data, and we 
evaluated the number of days that are required before all NPIs can be 
relaxed in relationship to: (i) risk tolerance; (ii)
maximum vaccination uptake; and (iii) daily vaccination rate. We uncover important aspects related to timing and coordination of the NPIs in the various regions based on traveling patterns.

\section*{Results}
\label{sec:results}

\begin{figure}[t!]
\begin{flushleft}
\small State-wide vaccination rate = $15,000$ vax/day:
\end{flushleft}
\vspace{-.4cm}
\subfigure[]{\includegraphics[width=.24\columnwidth]{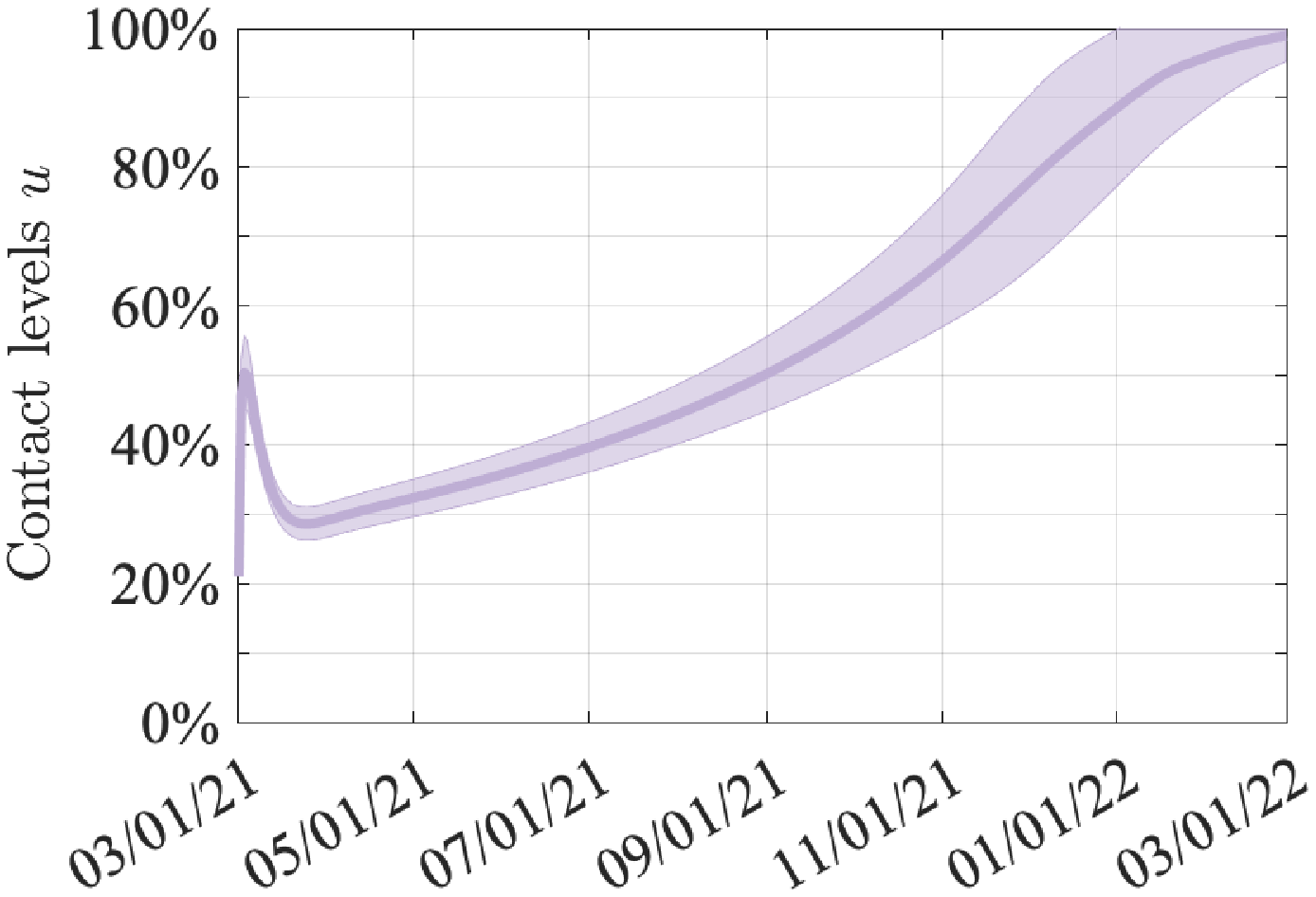}} 
\setcounter{subfigure}{2}%
\subfigure[]{\includegraphics[width=.24\columnwidth]{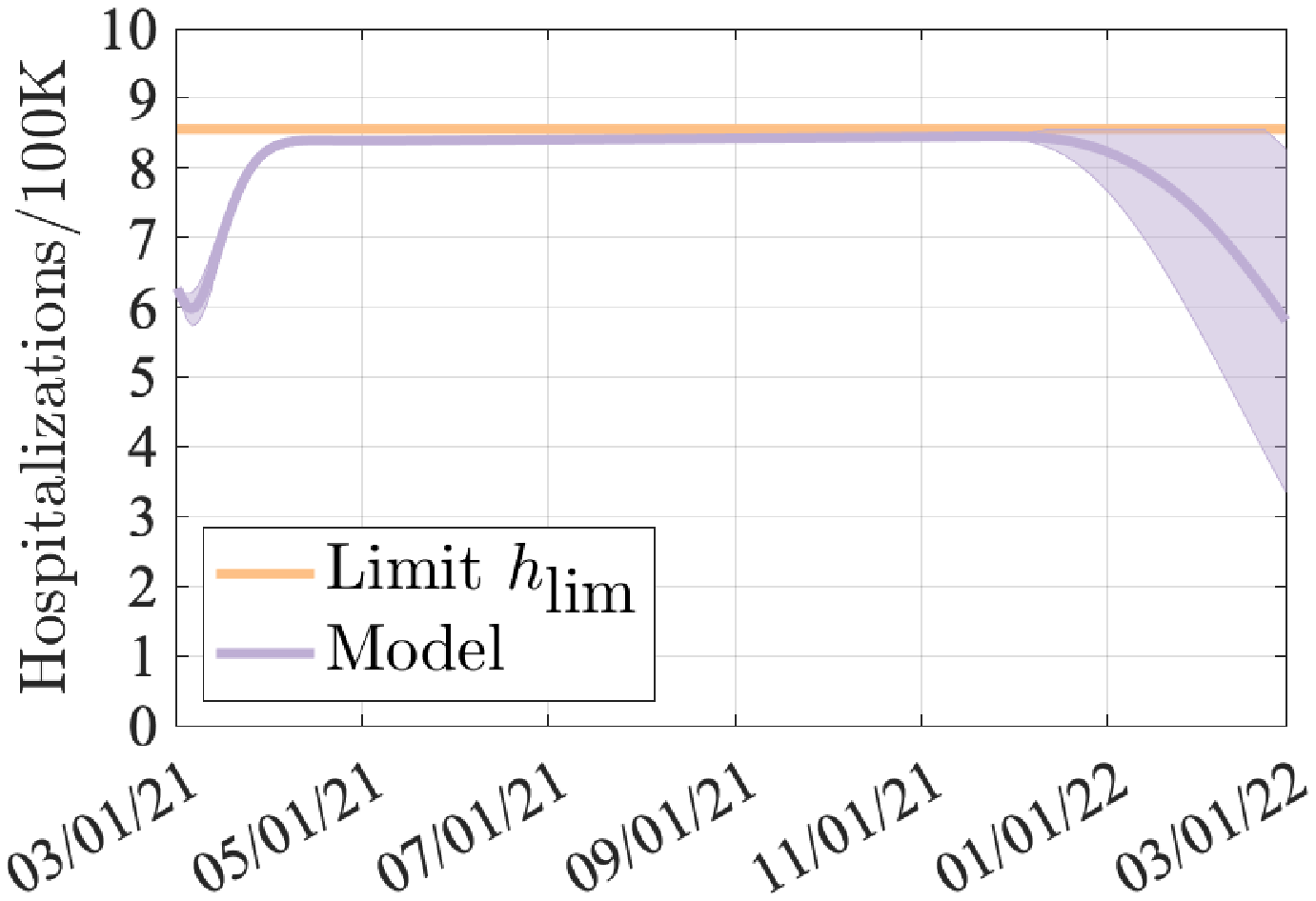}} 
\setcounter{subfigure}{4}%
\subfigure[]{\includegraphics[width=.24\columnwidth]{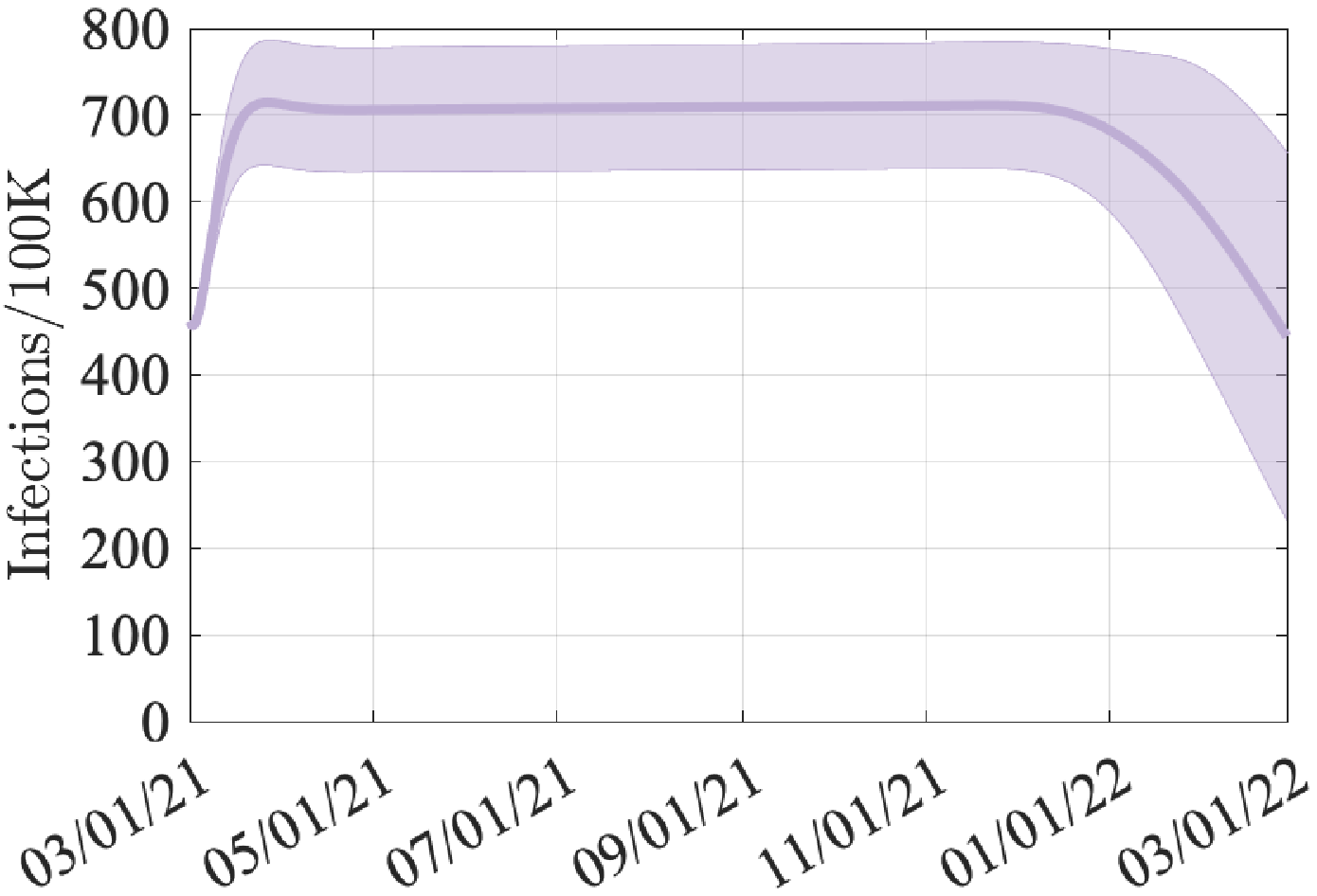}} 
\setcounter{subfigure}{6}%
\subfigure[]{\includegraphics[width=.24\columnwidth]{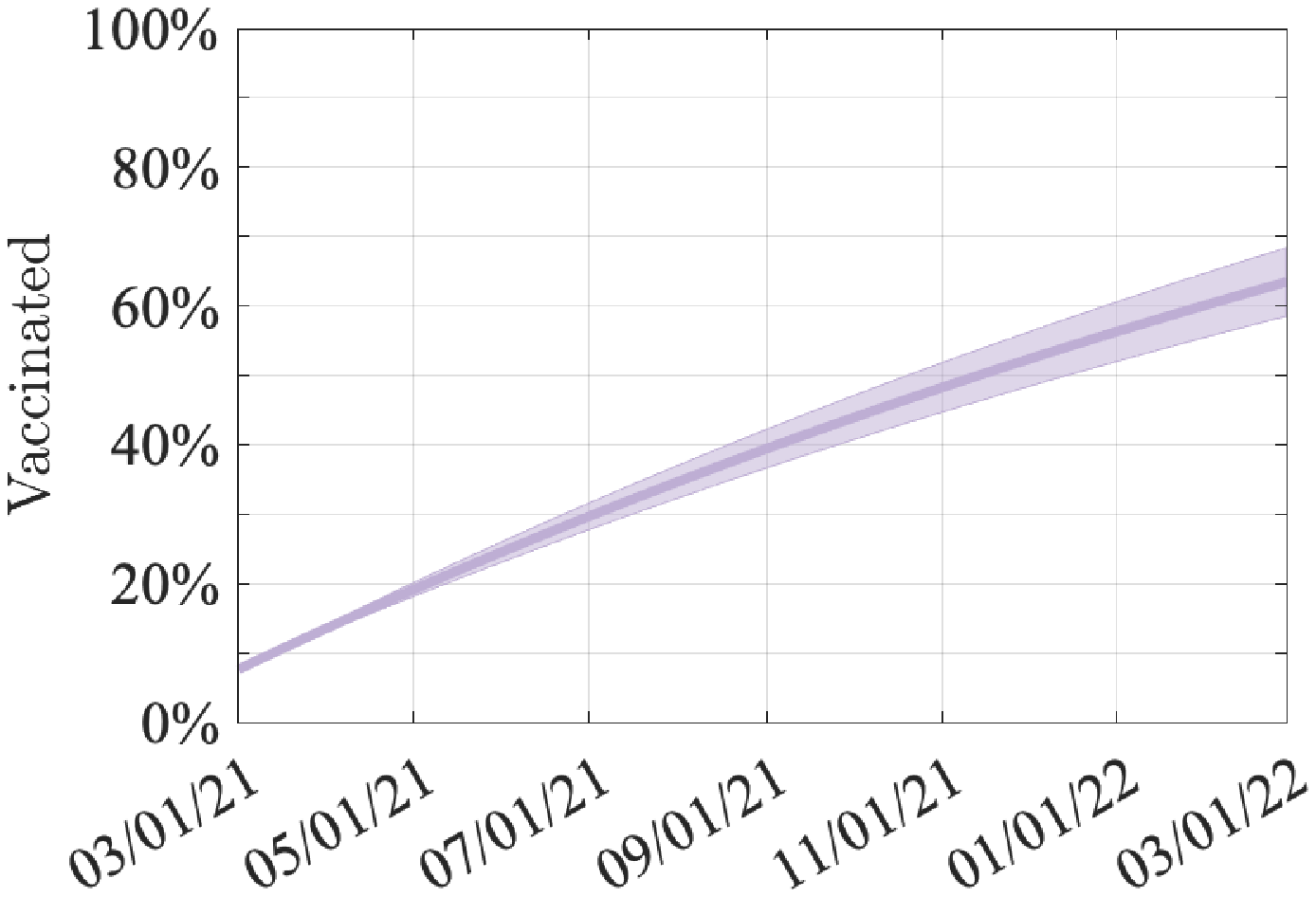}} 
\begin{flushleft}
\small  State-wide vaccination rate = $25,000$ vax/day:
\end{flushleft}
\vspace{-.4cm}
\setcounter{subfigure}{1}%
\centering \subfigure[]{\includegraphics[width=.24\columnwidth]{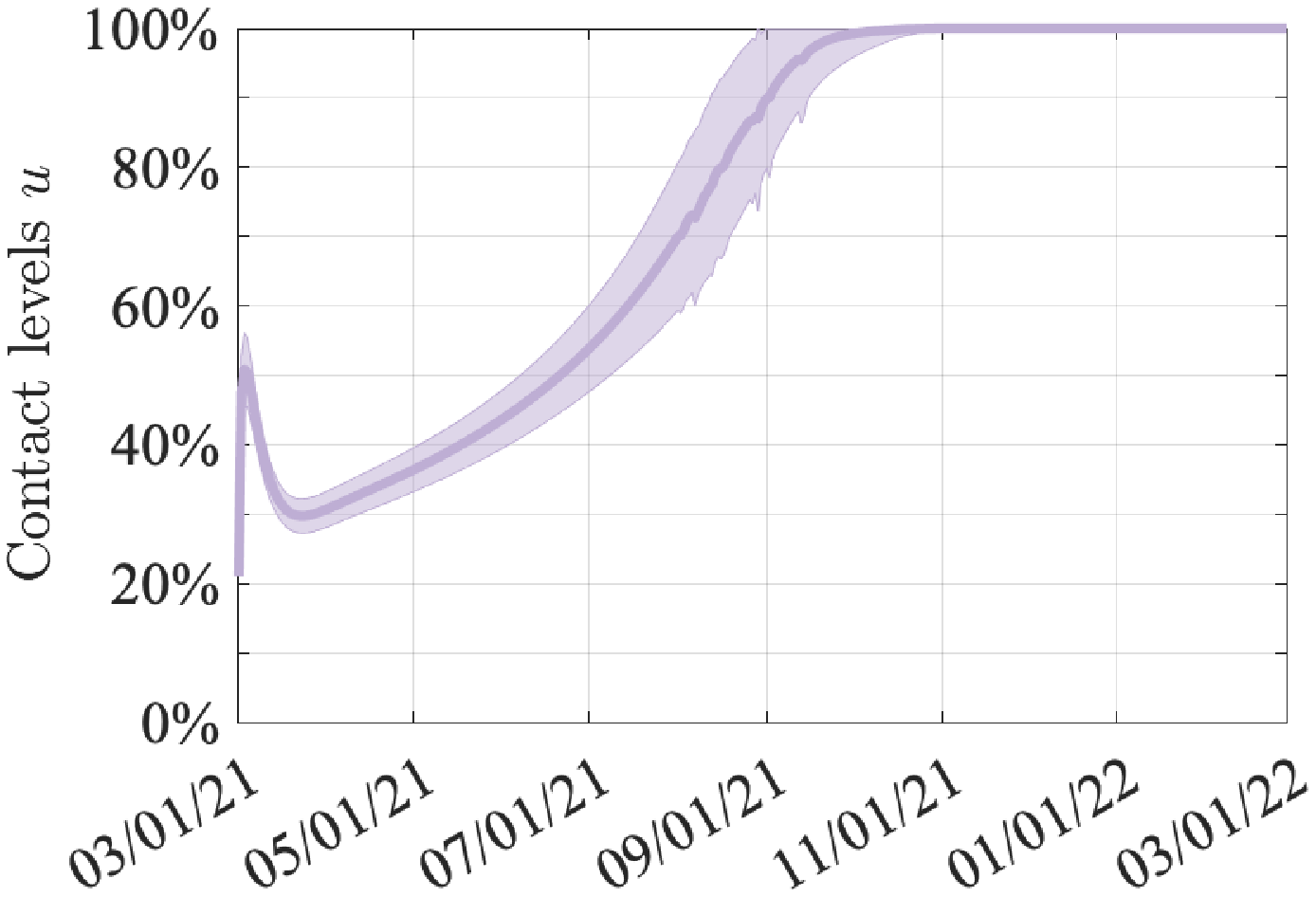}} 
\setcounter{subfigure}{3}%
\centering \subfigure[]{\includegraphics[width=.24\columnwidth]{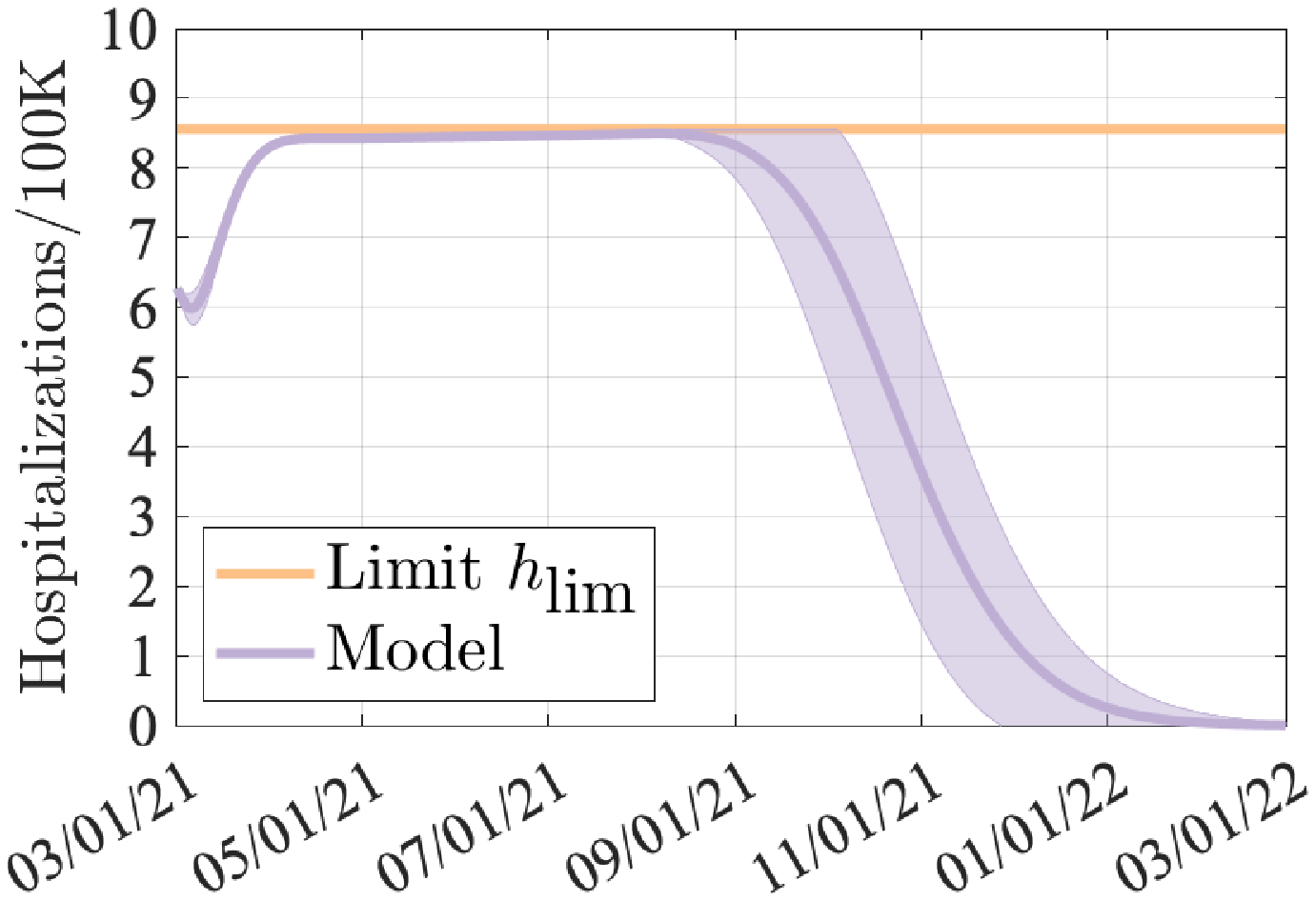}} 
\setcounter{subfigure}{5}%
\centering \subfigure[]{\includegraphics[width=.24\columnwidth]{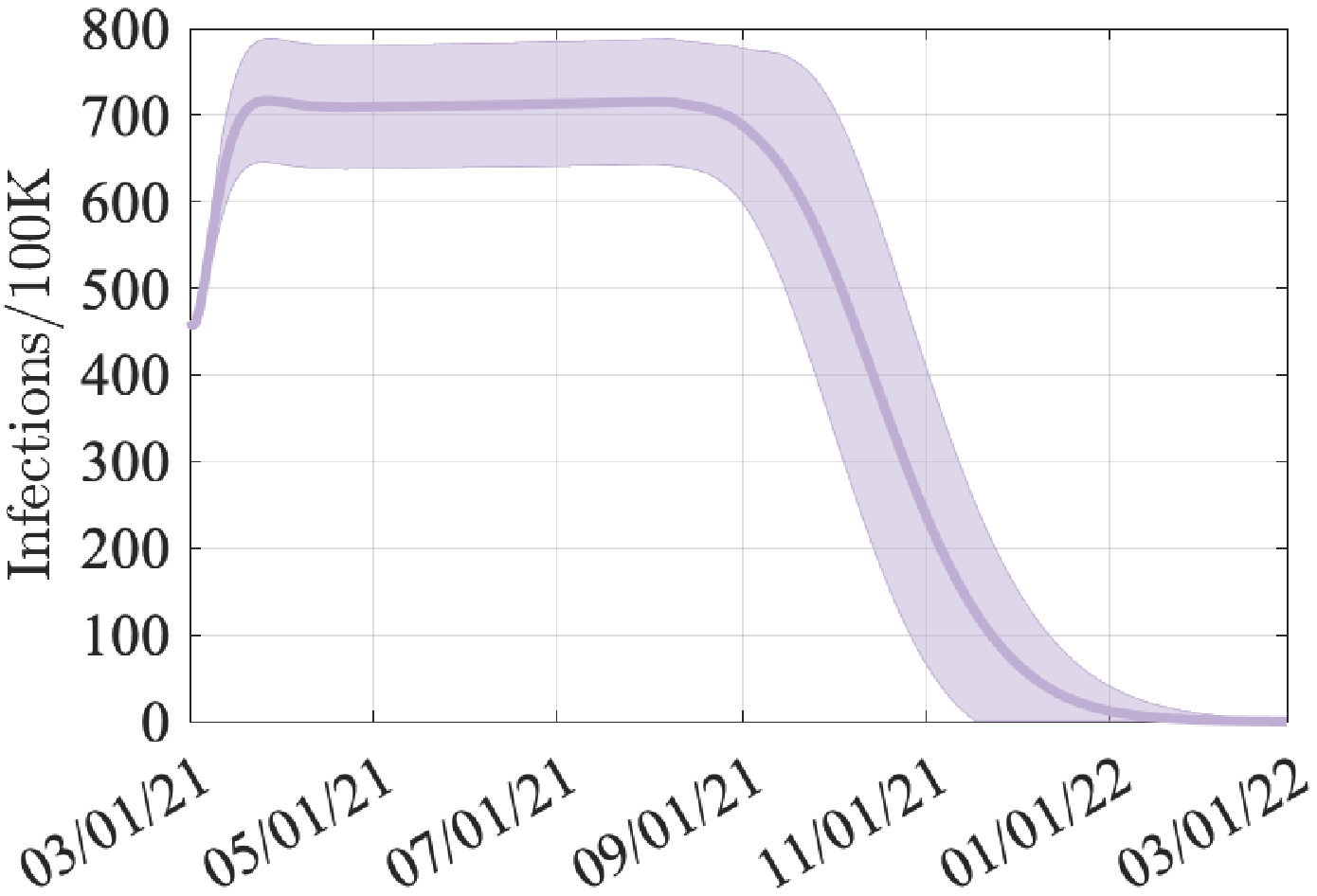}} 
\setcounter{subfigure}{7}%
\centering \subfigure[]{\includegraphics[width=.24\columnwidth]{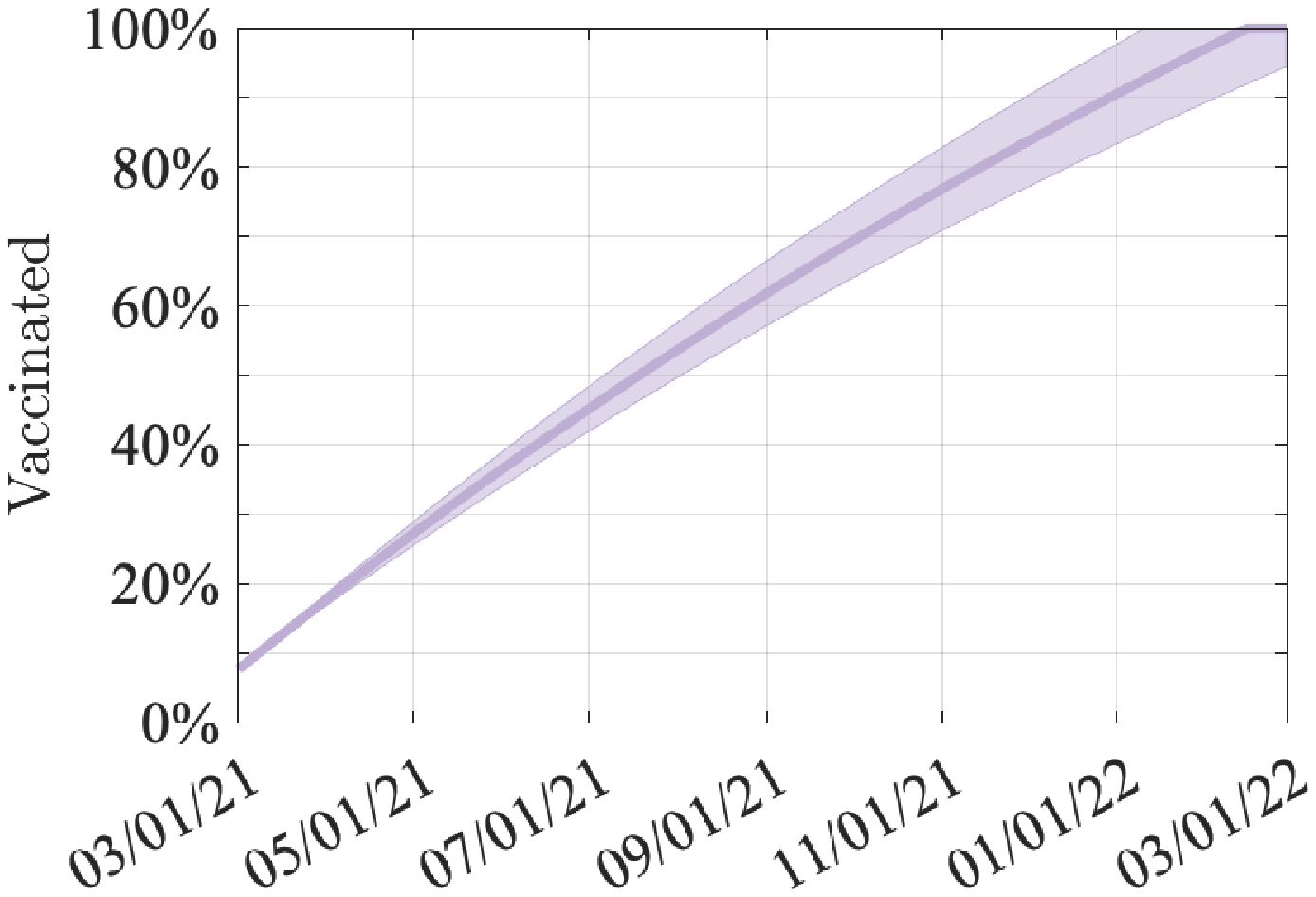}} 
\caption{\emph{\small{Model behavior when the feedback law is 
designed to simultaneously maximize contact levels and 
maintain hospitalizations below the threshold $\sbs{h}{lim}$.
(a)-(b) Level of transmission-relevant contacts with respect to 
pre-pandemic behavior, as selected by the feedback controller. 
All simulations are conducted by using a single-compartment model that 
is fitted using data from the state of Colorado,  USA (see  Methods). 
Results are averaged over $10,000$ simulations with parameters sampled 
using a Latin Hypercube technique within $15\%$ of their nominal values.
Continuous line shows mean of the trajectory and shaded area show 
$99.73\%$ confidence intervals. This figure shows an ideal situation 
where vaccination uptake can reach a level of $100\%$.}}}
\label{fig:trajectories_singleRegion_vax}
\end{figure}

\paragraph{Controlling Hospitalizations During Pandemic Outbreaks.} 
We aim to address the following question: ``What is the 
least-restrictive level of NPIs that guarantees that the number of 
hospitalized  individuals on each day does not exceed a pre-specified 
limit and simultaneously accounts for the economic implications of 
the NPIs?''. 
The hospitalization limit models the level of risk-tolerance in a 
population or the level of stress tolerated by the healthcare system in 
a given region.  
We answer this question by formulating a constrained optimization 
problem, which uses a compartmental model of epidemic transmission  
to predict the epidemic state (see Methods). The intensity of NPIs is
represented by a parameter $u \in[0,1]$ that describes the level of 
permitted transmission-relevant contacts, where $u=0$ corresponds to 
zero contacts (full lockdown) and $u=1$ corresponds to pre-pandemic 
contact levels (zero NPIs) (see 
\figurename~\ref{fig:trajectories_singleRegion_vax}(a)--(b)).
An optimization problem is used to derive a feedback law that uses the 
instantaneous epidemic state to systematically select a level of 
transmission-relevant contacts that balances between the economic impact
of the imposed restrictions and the number of infectious individuals, 
while simultaneously guaranteeing that the number of daily 
hospitalizations does not exceed the specified hospitalization 
limit (\figurename~\ref{fig:trajectories_singleRegion_vax}(c)--(d)), 
denoted by $\sbs{h}{lim}$. Here, $\sbs{h}{lim}$ models the maximum 
allowable number of hospitalized individuals on each day.
The parameters of our model are fitted 
to official data from the state of Colorado, USA  
\cite{AB-EC-DG-etal-b:21,AB-EC-DG-etal:21}.

\paragraph{When Can We Safely Return to Normal?} 
To address this question, we used the feedback optimization framework 
to determine the highest allowable level of contacts (on each day) that 
ensures that the daily number of hospitalized individuals does not 
exceed (and remains close to) the pre-specified limit $\sbs{h}{lim}$. 
In a single simulation (
\figurename~\ref{fig:trajectories_singleRegion_vax}) it can be seen that, as 
the fraction of vaccinated individuals in the population increases, the 
allowable contact levels selected by the feedback law can also gradually
increase, while ensuring that the number of hospitalizations remains 
below the pre-defined limit; the infection rate is similarly 
constrained.  
Accordingly, our framework allows us to characterize a lower bound on
the number of days required before a full return to normality can safely
occur.

\begin{figure}[t!]
\centering
\centering \subfigure[]{\includegraphics[width=.34\columnwidth]{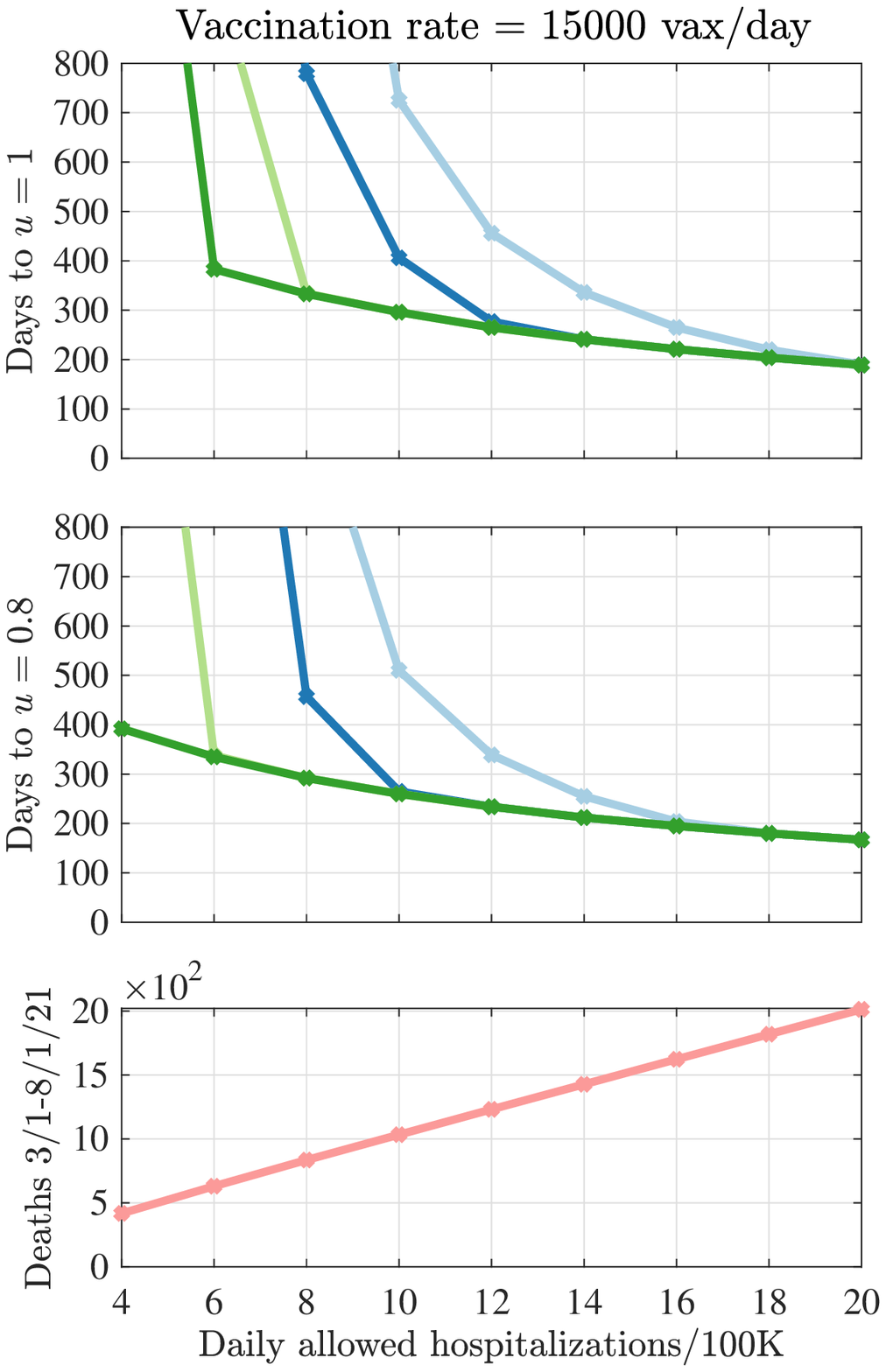}}%
\hfill
\centering \subfigure[]{\includegraphics[width=.3\columnwidth]{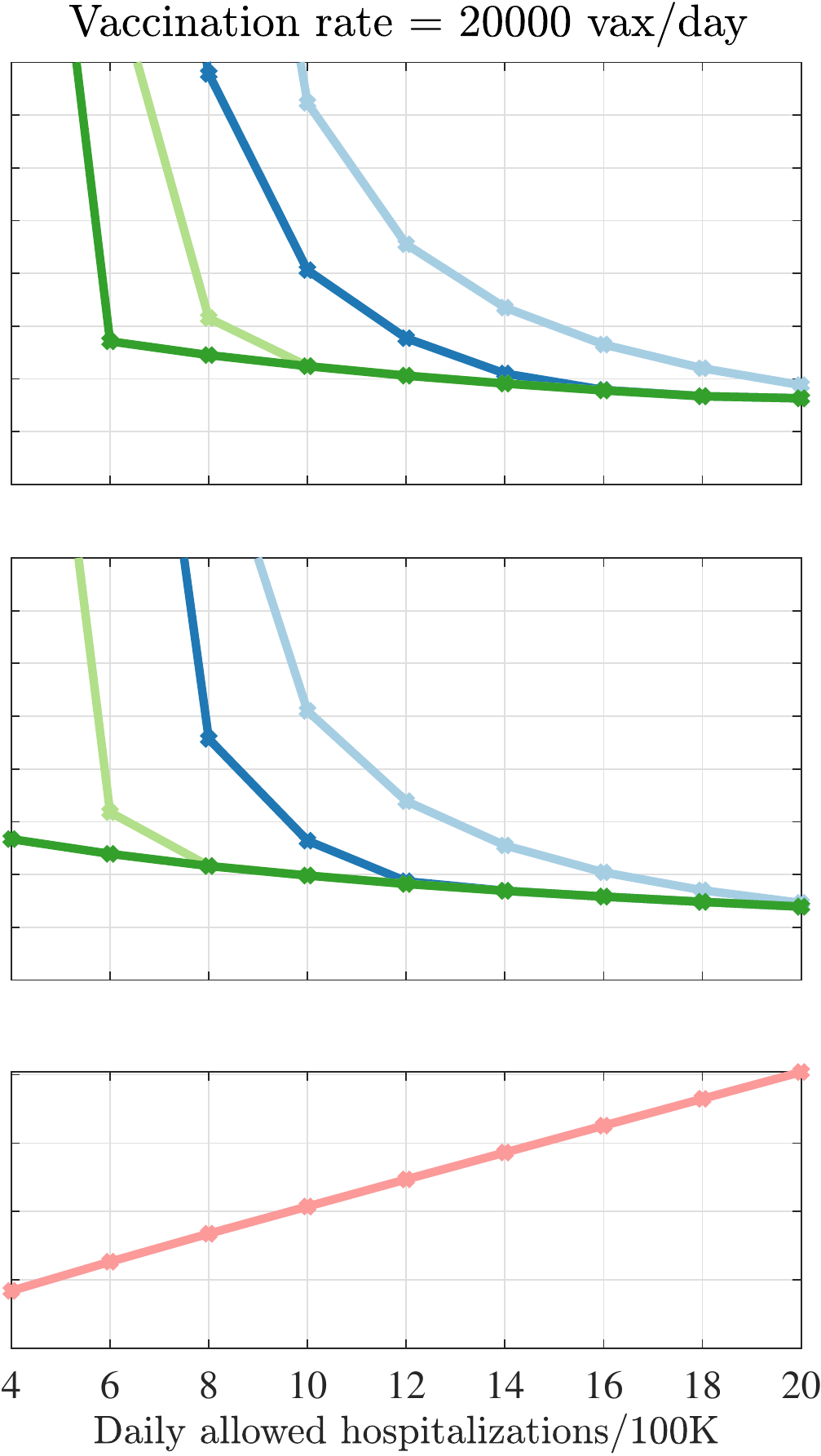}}%
\hfill
\centering \subfigure[]{\includegraphics[width=.3\columnwidth]{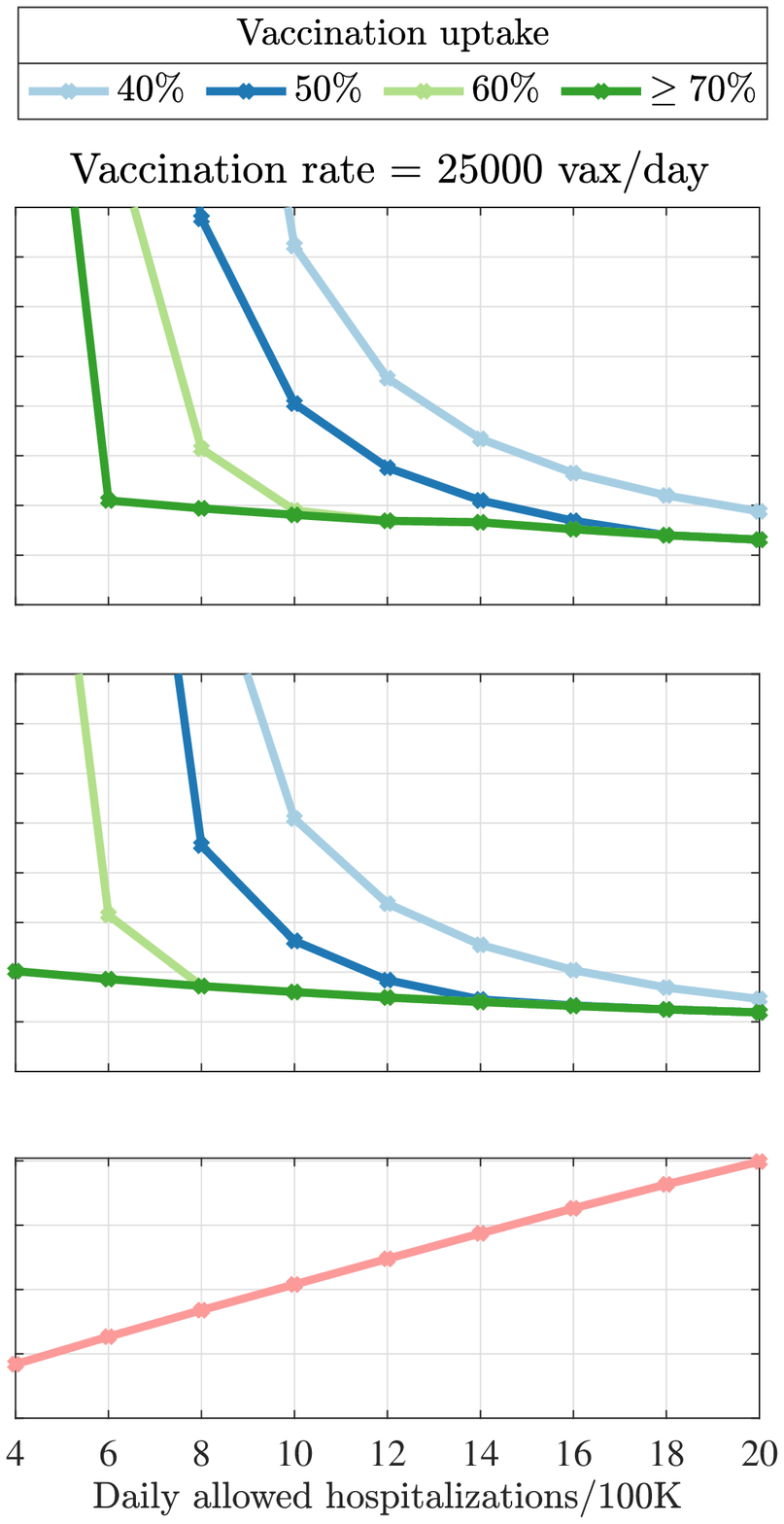}}\\
\caption{\emph{\small{
Number of days needed before a return to normal can be implemented 
without exceeding predefined hospitalization limits. The number of days is counted 
beginning 3/1/21.
(Top row:) Number of days to  $u=1$. 
(Center row:) Number of days to  $u=0.8$.  
(Bottom row:) Estimated number of deaths between March 1, 2021, and 
August 1, 2021. All simulations are conducted by using a 
single-compartment model fitted using data from the state of 
Colorado,  USA (see  Methods). Any vaccination uptake of $70\%$ or 
larger yields an identical (dark green) curve.}}}
\label{fig:daysTou=1}
\end{figure}

Return to normality refers to a condition where all NPIs can be 
repealed and the societal behavior can return to pre-pandemic contact 
levels  (i.e. $u=1$).
We found that the parameter that most consistently impacted the number 
of days to normality is the vaccination uptake in the population (see  
\figurename~\ref{fig:daysTou=1} where the number of days to normality is counted 
beginning March 1 st, 2021). This finding follows 
from three main observations. 
First, our results show that any vaccination uptake of $50\%$ or less will 
require more than $2$ years (730 days) to return to normal behavior for any of the examined vaccination rates if we expect to keep the 
number of hospitalizations below 8 individuals/day per 100K inhabitants.
Second, the number of days to normality reduces by a factor of at least $2.5$ 
as the vaccination uptake is increased from $50\%$ to $60 \%$
(for instance, with vaccination rate of $25,000$ vax/day and hospitalization 
limit $\sbs{h}{lim}=8$, the number of days to $u=1$ decreases from $777$ to  
$314$ as the vaccination uptake is varied from $50\%$ to $60 \%$).
Third, our results suggest that vaccination uptakes larger that $70\%$
will not decrease the time before all NPIs can be safely lifted. This is likely a
factor of the current level of infection-derived immunity in the population, such
that $70\%$ vaccine uptake is sufficient to reach a threshold where infections 
decline, regardless of contact levels. 
A second parameter that consistently affects the time to normality is 
the level of risk tolerance $\sbs{h}{lim}$, i.e., the number of hospitalizations 
tolerated during the outbreak. 
Precisely, allowing more hospitalizations to occur (increasing $\sbs{h}{lim}$) 
reduces the number of days to normality. Decreased risk tolerance leads to a 
lower infection rate and decreases the rate of naturally-acquired immunity, thus 
increasing the time required for a  return to normality.
For instance, with a vaccination rate of $15,000$ individuals/day, when vaccination uptake is at least $70\%$, the number of 
days to $u=1$ reduces from over $1$ year ($383$ days) to about $6$ months
($189$ days) when the hospitalization limit is increased from $6$ to $20$ 
individuals/day.
Unfortunately, although a higher hospitalization limit reduces the 
time to normality, it results in a higher number of deaths 
(\figurename~\ref{fig:daysTou=1} bottom panel). Our results also suggest that 
the vaccination uptake does not affect the number of deaths in the considered 
time interval (this fact emerges because the number of deaths is counted until 
August 1, which occurs before the vaccination uptake threshold is reached).

Our results suggest that in order to return to normal behavior (pre-pandemic contact rates) on 10/01/21 (153 days beginning 3/1/21), there is likely to be stress on the healthcare system (i.e. the number of hospitalizations will exceed the pre-defined limit), the level of which will depend on vaccine uptake. If vaccine uptake is as low as $40\%$, we may not remain below $20$ hospitalizations/day per 100K inhabitants. However, with vaccine uptake of $60\%$ we could return to normal and keep hospitalizations below $16$ hospitalizations/day per 100K inhabitants.

In practice, we expect that contact behaviors will vary widely across 
individuals and that, even when all government policies will be repealed, many 
individuals will likely continue to practice NPIs, such as mask wearing and 
social distancing.
Thus, although $u=1$ may be an unlikely scenario in the near future, individuals 
may quickly resume a behavior of ``almost normality'', where contacts are 
restored to $80\%$ of pre-pandemic levels (i.e., $u=0.8$).

Following that, in order to safely return to $u=0.8$ on 10/01/21, if vaccine uptake is as low as $50\%$, we may not remain below $14$ hospitalizations/day per 100K inhabitants. However, with vaccine uptake of $70\%$ we could return to normal and keep hospitalizations below $12$ hospitalizations/day per 100K inhabitants. 
(see \figurename~\ref{fig:daysTou=1} center row).

\paragraph{Effects of Regional Heterogeneities and Mobility Patterns.}

\begin{figure}[t!]
\centering 
\includegraphics[width=1.0\columnwidth]{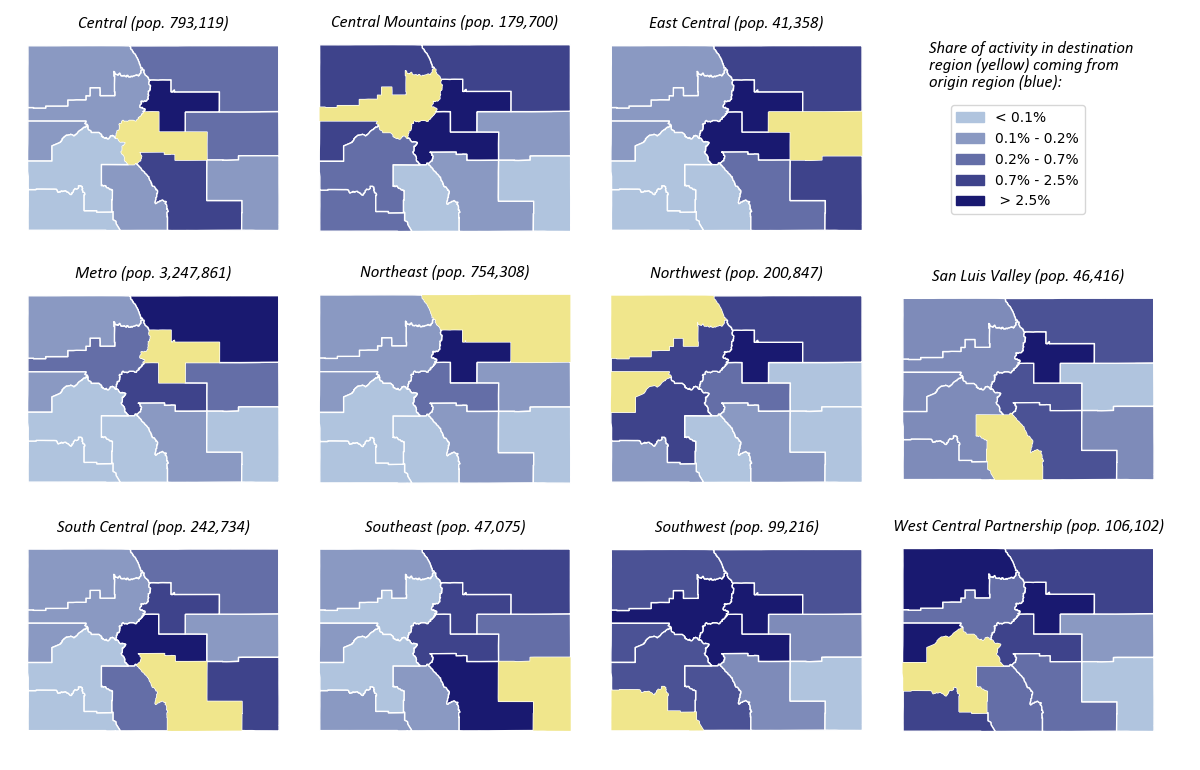}
\caption{\emph{ \small{Regional connectivity patterns between the 11 Local 
Public Health Agency regions in Colorado, USA. Each panel illustrates the 
intensity of contacts between residents of the yellow region and individuals
traveling from the blue regions. 
Total travel volume is averaged over the time period 1/1/20--12/31/20.
Data obtained from Safegraph.}}}
\label{fig:network_regions_graph}
\end{figure}

Heterogeneities in policies and  
behavioral responses to an ongoing epidemic between regions motivate the use of 
higher-resolution models and control techniques that can adequately 
capture this diversity. In this context, a crucial open 
question concerns how local authorities can identify region-dependent 
levels of NPIs that guarantee that local hospitalization limits are 
met, and to what extent inter-regional decision-making can be 
coordinated to achieve this objective.
To address these questions, we generalized the transmission model and 
feedback control framework to a network setting (see Methods), and we 
used publicly-available mobility data from cell-phone usage to estimate
inter-regional couplings (\figurename~\ref{fig:network_regions_graph} 
illustrates regional connectivity patterns for the state of Colorado, 
USA. Note that the Metro region comprises the majority of the state's 
population and contributes a large fraction of activity statewide.
The model is organized into eleven regions, each describing a Local 
Public Health Agency (LPHA) region in Colorado, USA 
\cite{AB-EC-DG-etal:21}, and transmission levels are fitted to 
regional hospitalization data from the period 1/1/21--2/28/21. 

\begin{figure}[th!]
\subfigure[]{\includegraphics[width=.24\columnwidth]{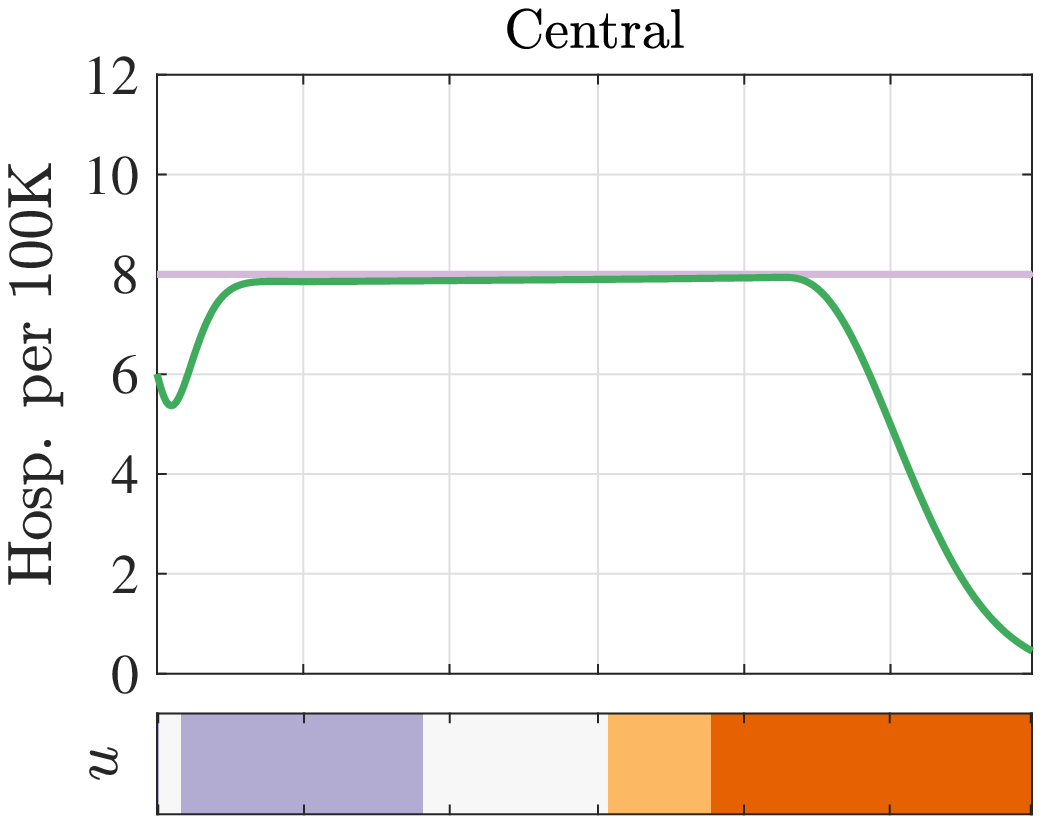}}%
\hspace{.01cm}
\subfigure[]{\includegraphics[width=.24\columnwidth]{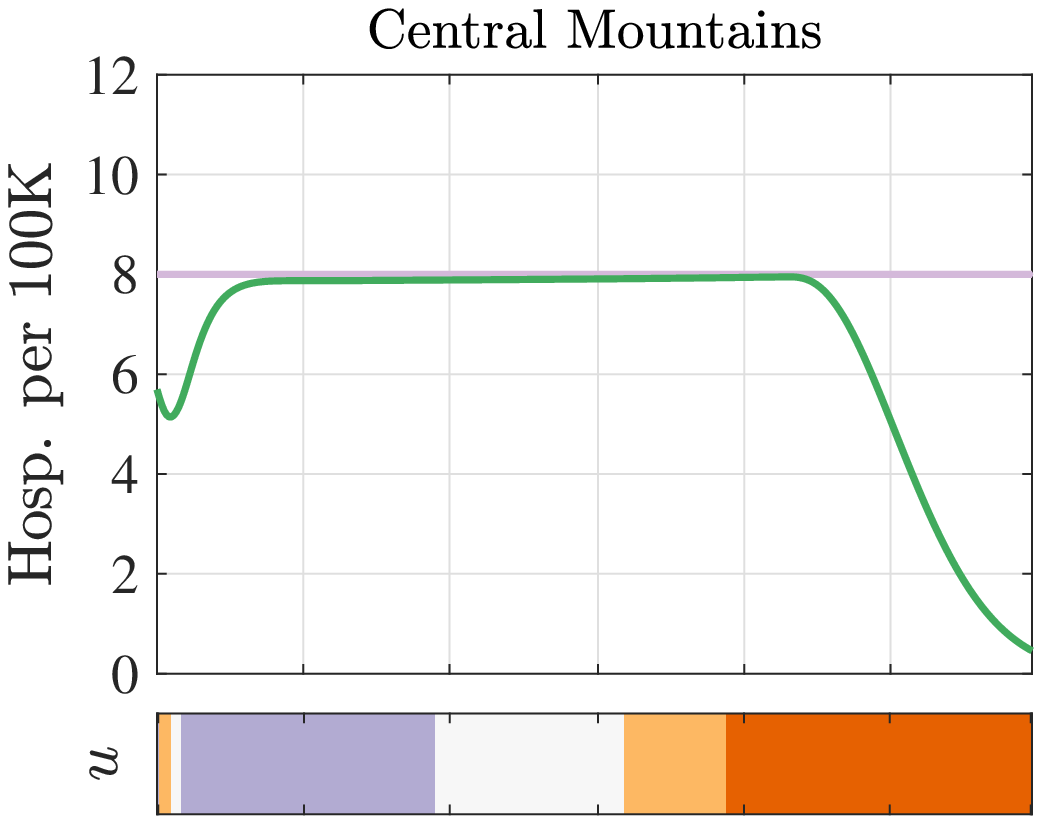}}%
\hspace{.01cm}
\subfigure[]{\includegraphics[width=.24\columnwidth]{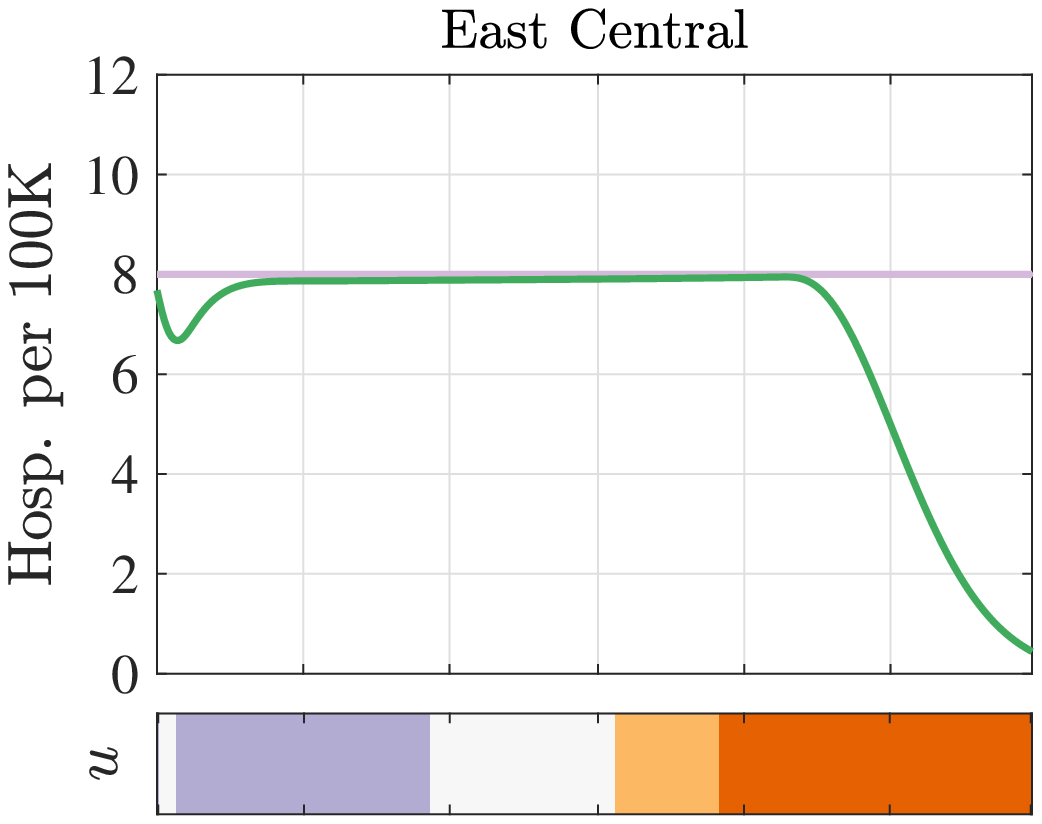}}%
\hspace{.1cm}
\subfigure[]{\includegraphics[width=.14\columnwidth]{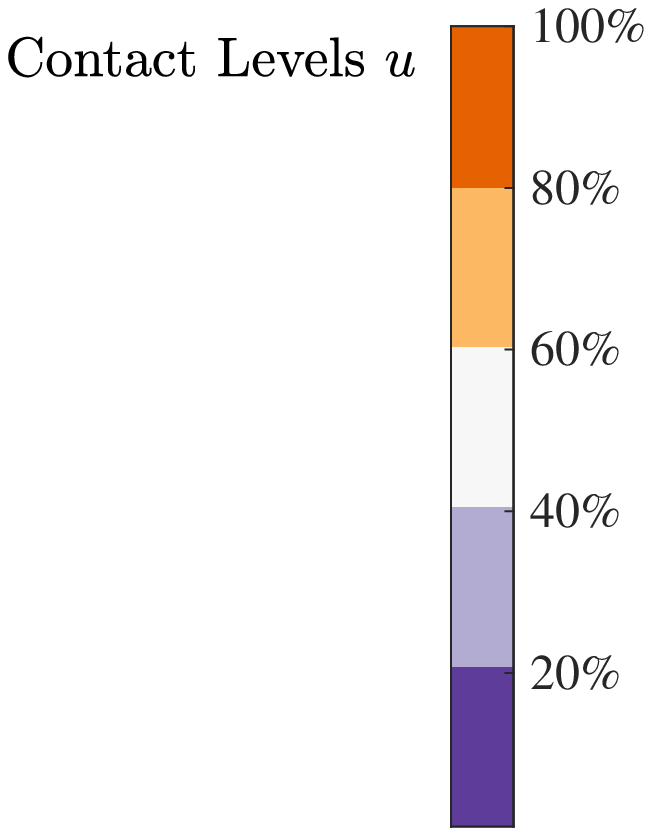}}\\
\subfigure[]{\includegraphics[width=.24\columnwidth]{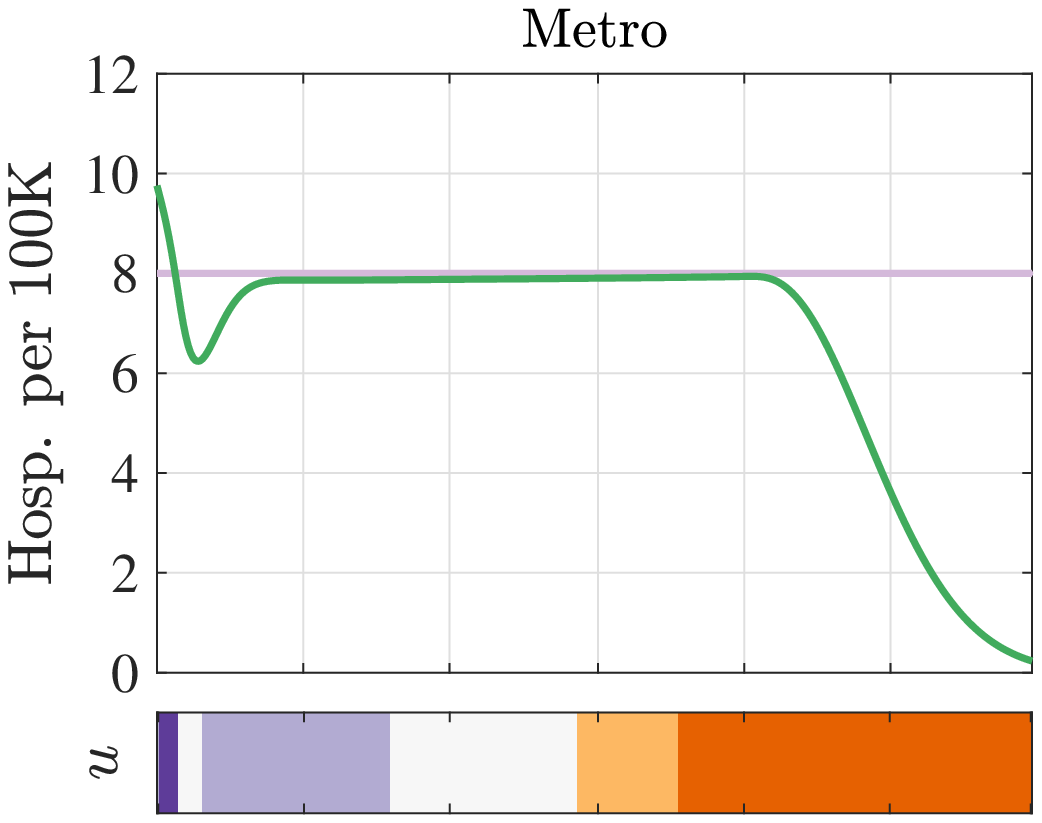}}%
\hspace{.01cm}
\subfigure[]{\includegraphics[width=.24\columnwidth]{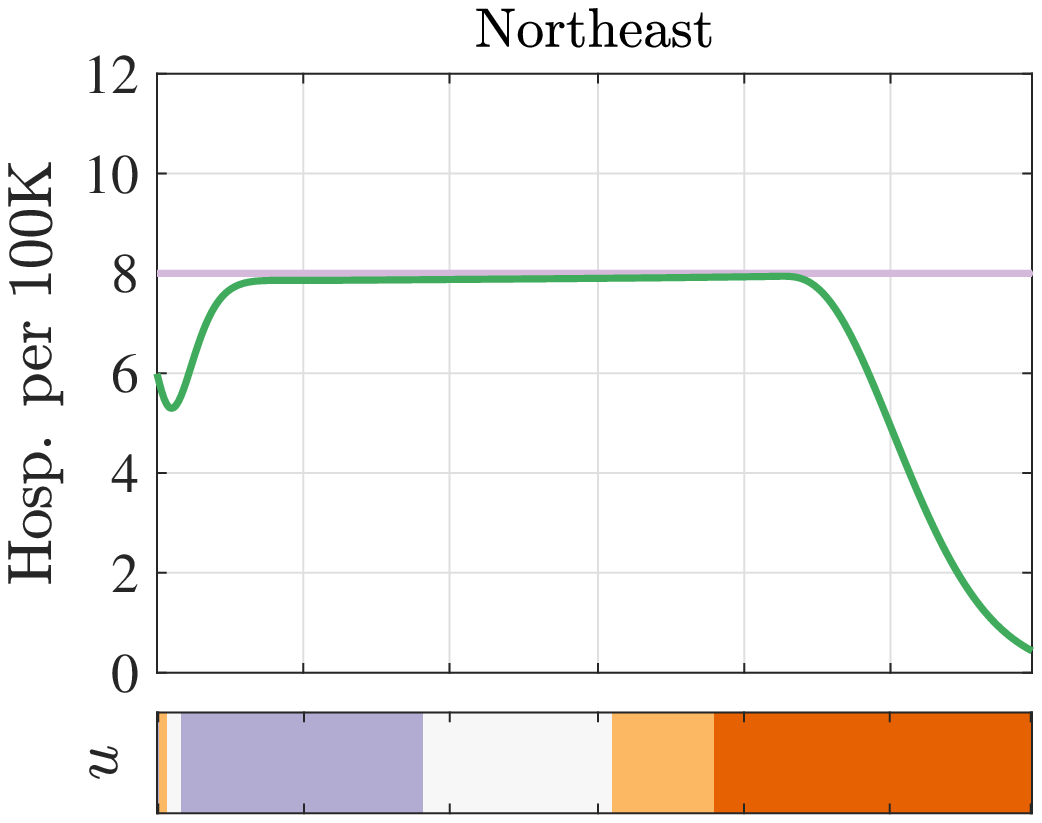}}%
\hspace{.01cm}
\subfigure[]{\includegraphics[width=.24\columnwidth]{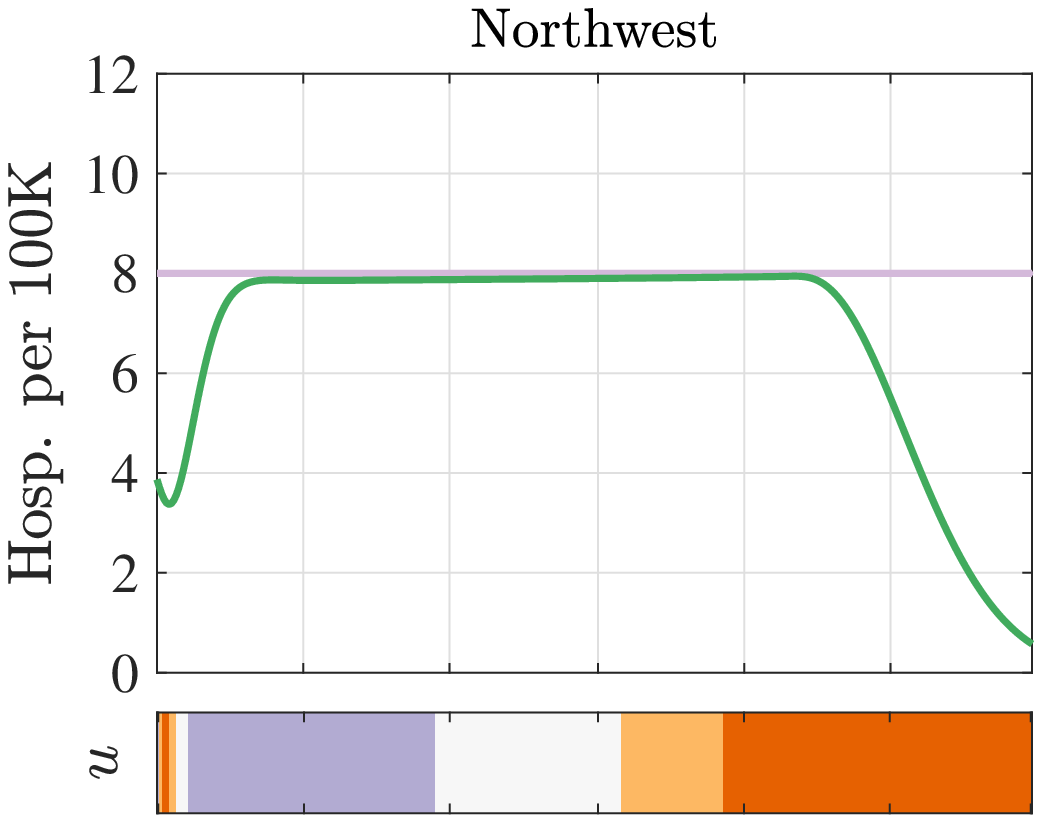}}%
\hspace{.01cm}
\subfigure[]{\includegraphics[width=.24\columnwidth]{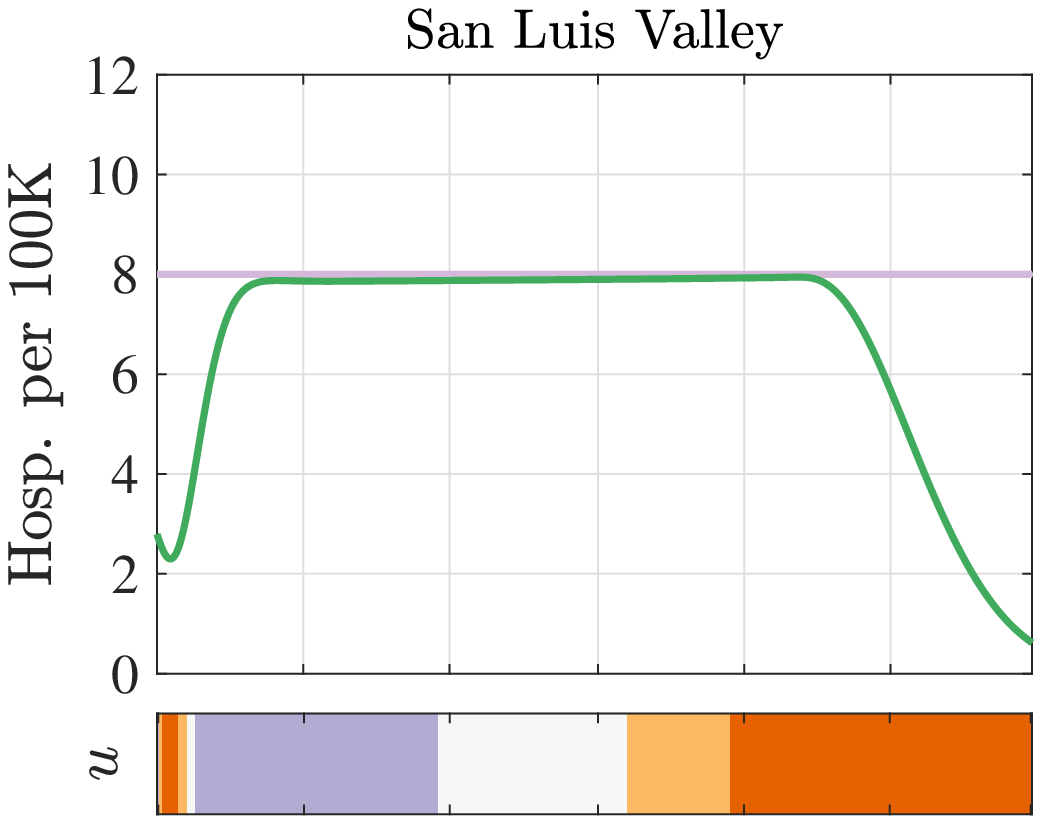}}\\
\subfigure[]{\includegraphics[width=.24\columnwidth]{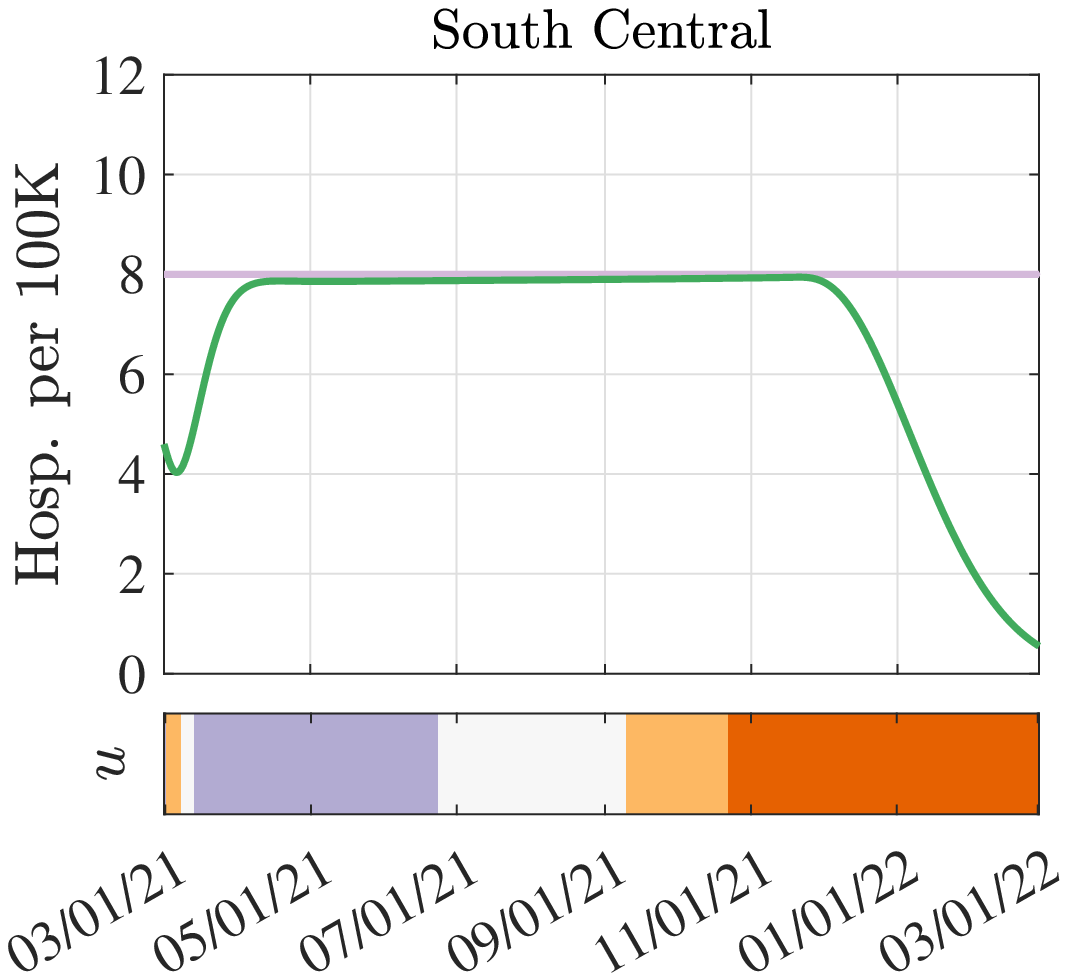}}%
\hspace{.01cm}
\subfigure[]{\includegraphics[width=.24\columnwidth]{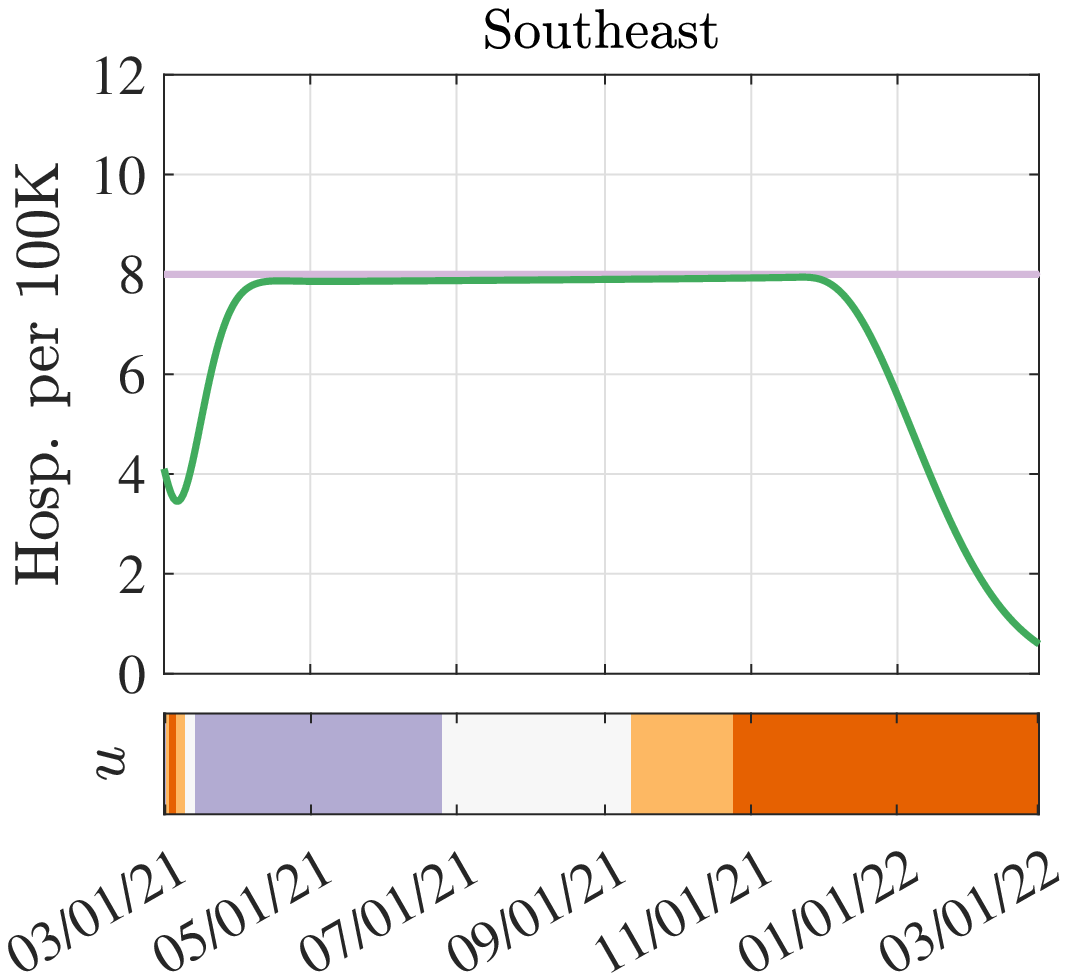}}%
\hspace{.01cm}
\subfigure[]{\includegraphics[width=.24\columnwidth]{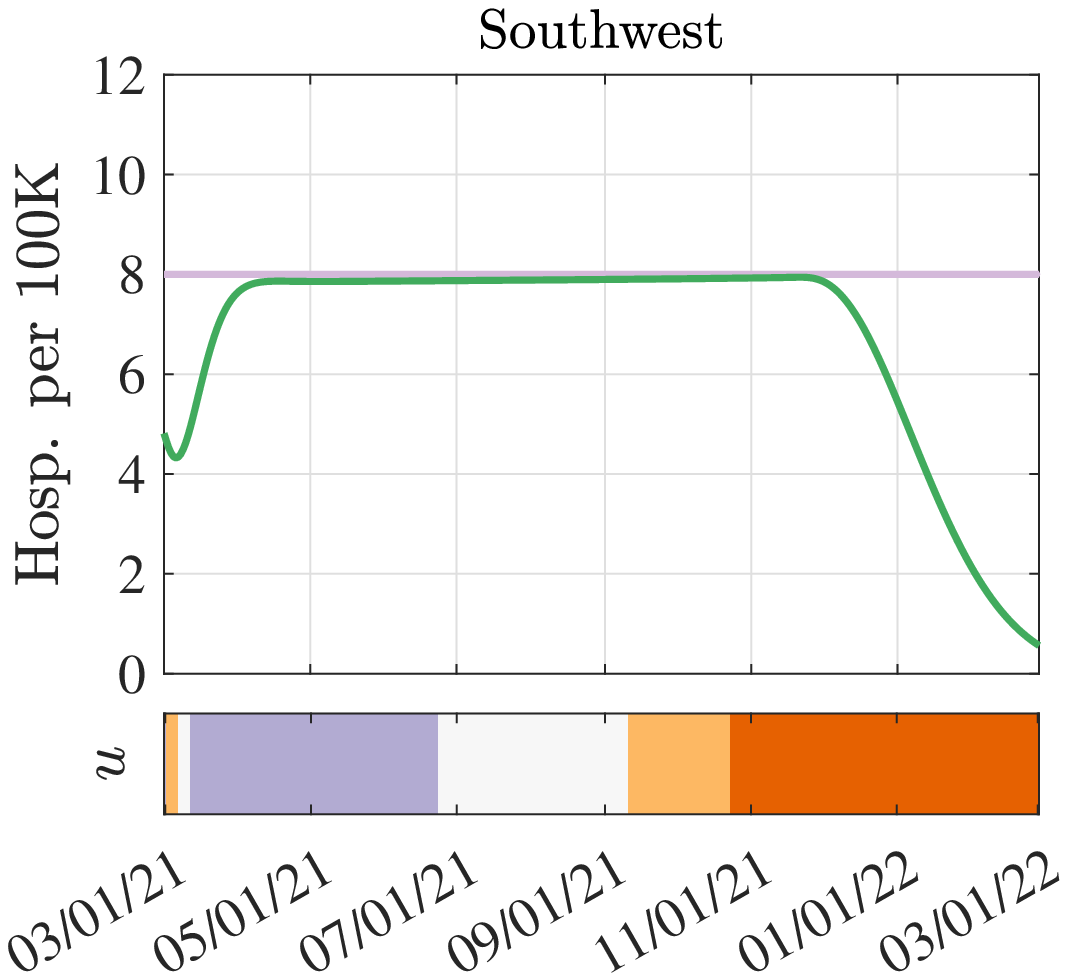}}%
\hspace{.01cm}
\subfigure[]{\includegraphics[width=.24\columnwidth]{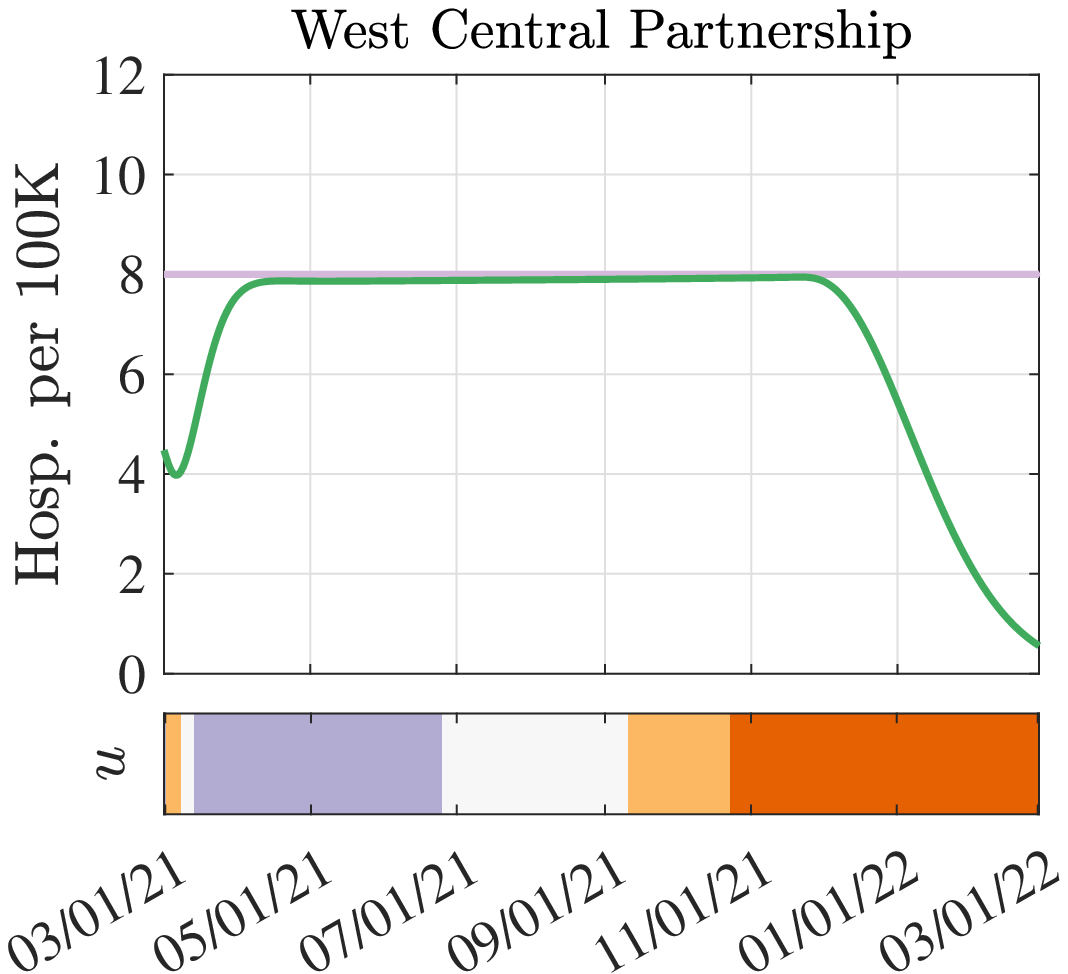}}%
\caption{\emph{\small{Hospitalizations and controller level over time when a 
group of regional controllers are used to guarantee that pre-specified 
region-dependent hospitalization limits are not violated. 
Each of the 11 panels shows the evolution in a different LPHA region 
(see \figurename~\ref{fig:network_regions_graph} for an illustration of 
the connectivity graph).
Simulation are conducted with a state-wide vaccination rate of 
$20,000$ vax/day, to a maximum vaccine update of $70\%$.
Solid green lines illustrate the evolution of the hospitalized state, 
light purple lines show the pre-specified hospitalization limit. 
Heat maps illustrate the required level of NPIs $u_i$, as determined by
the controller.}}}
\label{fig:regional1}
\end{figure}

Local levels of NPIs in each region $i$ are represented by a parameter
$u_i \in [0,1]$, describing the level of permitted 
transmission-relevant contacts (or NPIs) in the region.
By using the feedback-optimizing control law (see Methods), we derived 
region-dependent daily levels of  NPIs that guarantee that the number
of hospitalized individuals in each region does not exceed a
region-dependent hospitalization limit $\sbs{h}{lim,$i$}$.
\figurename~\ref{fig:regional1} illustrates the number of hospitalized
individuals in the time-interval 03/01/21 -- 03/01/22, together with 
the contact levels $u_i$, as selected by the controller.
By comparing the simulation outcomes across the various regions, our 
model and control methods suggest that heterogeneities among the regions can
be exploited by the feedback controller during the initial phase 
in order to regulate the number of hospitalizations near the specified 
threshold in all regions. 
After the initial phase, uniform levels of NPIs among the regions can be 
used. Together, our results indicate that regional heterogeneities can be 
used by the feedback controller, especially when the epidemic state varies 
widely across the regions. Note that control levels and  epidemic curves actually varied drastically between regions throughout the pandemic, indicating, among other heterogeneities, the difficulty in agreeing upon tolerable risk levels in a heterogeneous region.

\begin{figure}[th!]
\subfigure[]{\includegraphics[width=.24\columnwidth]{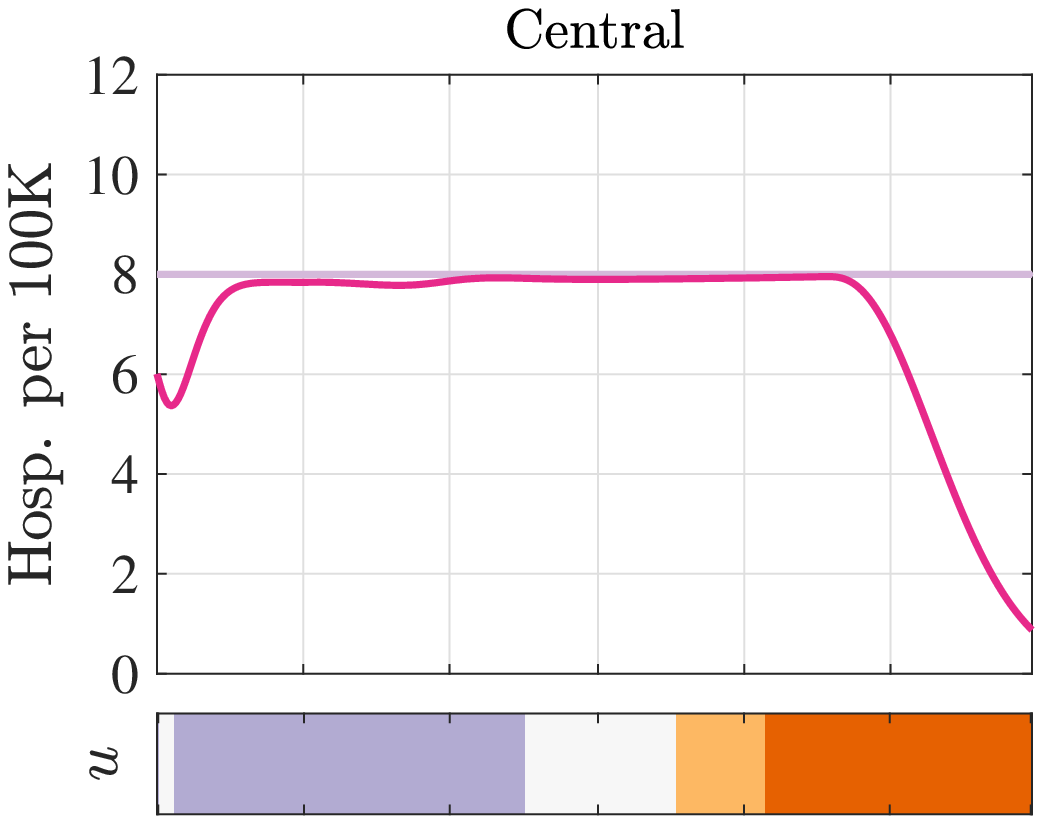}}%
\hspace{.01cm}
\subfigure[]{\includegraphics[width=.24\columnwidth]{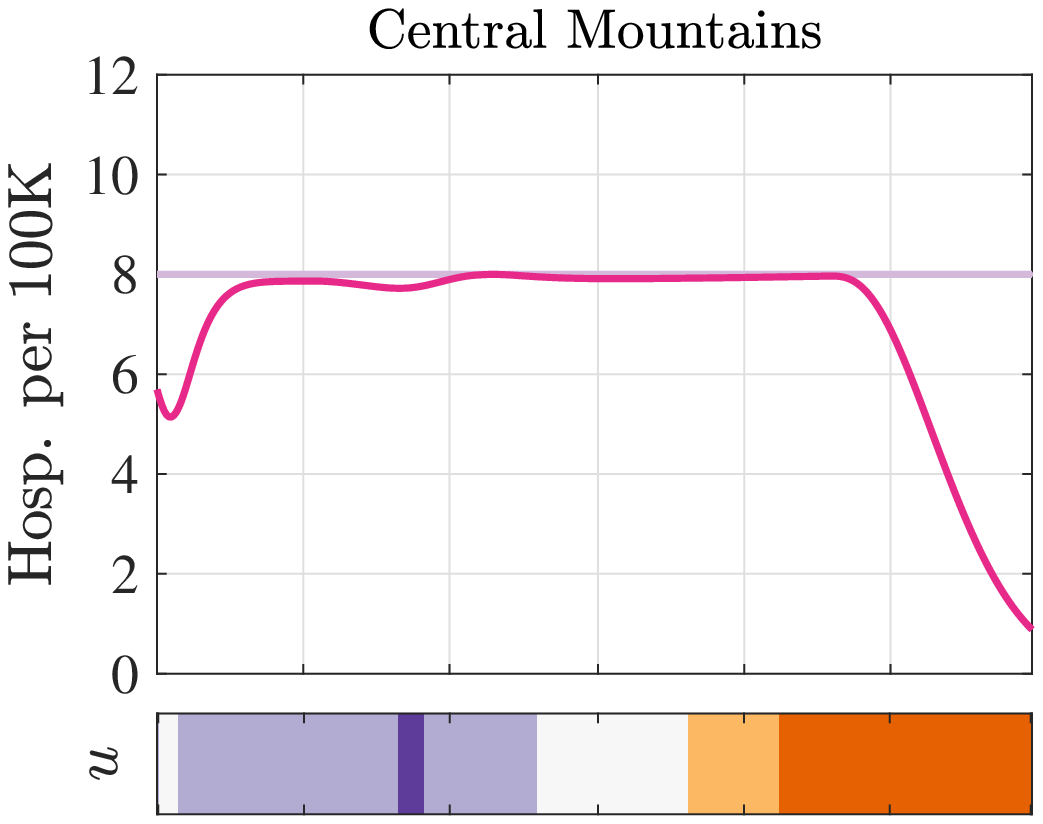}}%
\hspace{.01cm}
\subfigure[]{\includegraphics[width=.24\columnwidth]{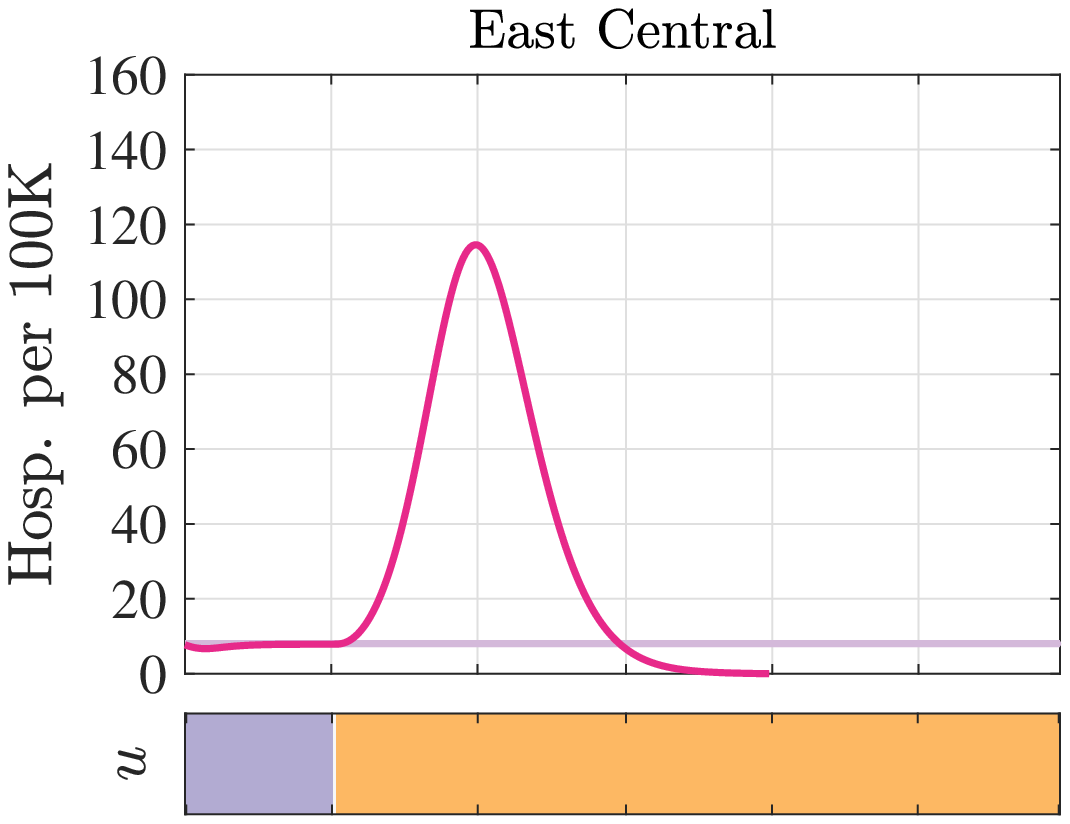}}%
\hspace{.01cm}
\subfigure[]{\includegraphics[width=.14\columnwidth]{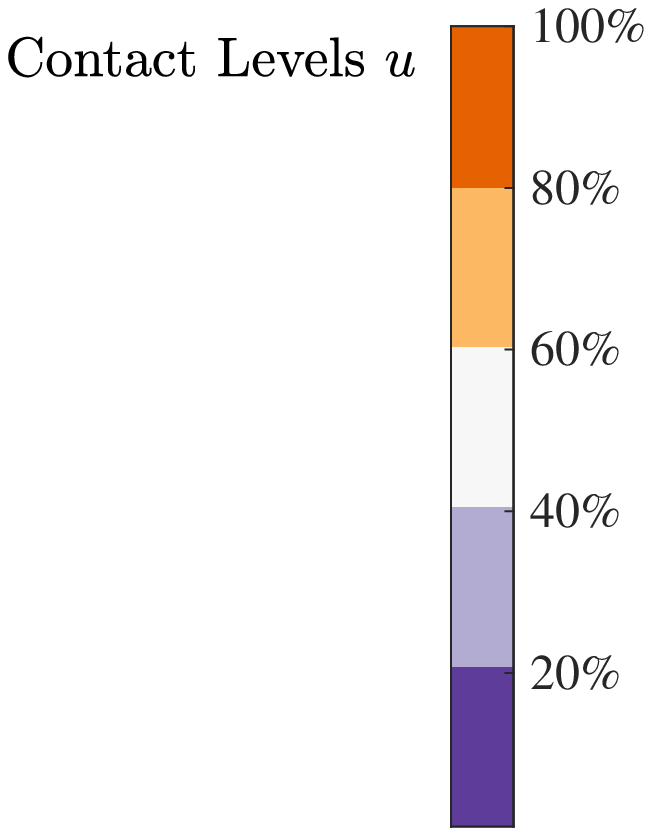}}\\
\subfigure[]{\includegraphics[width=.24\columnwidth]{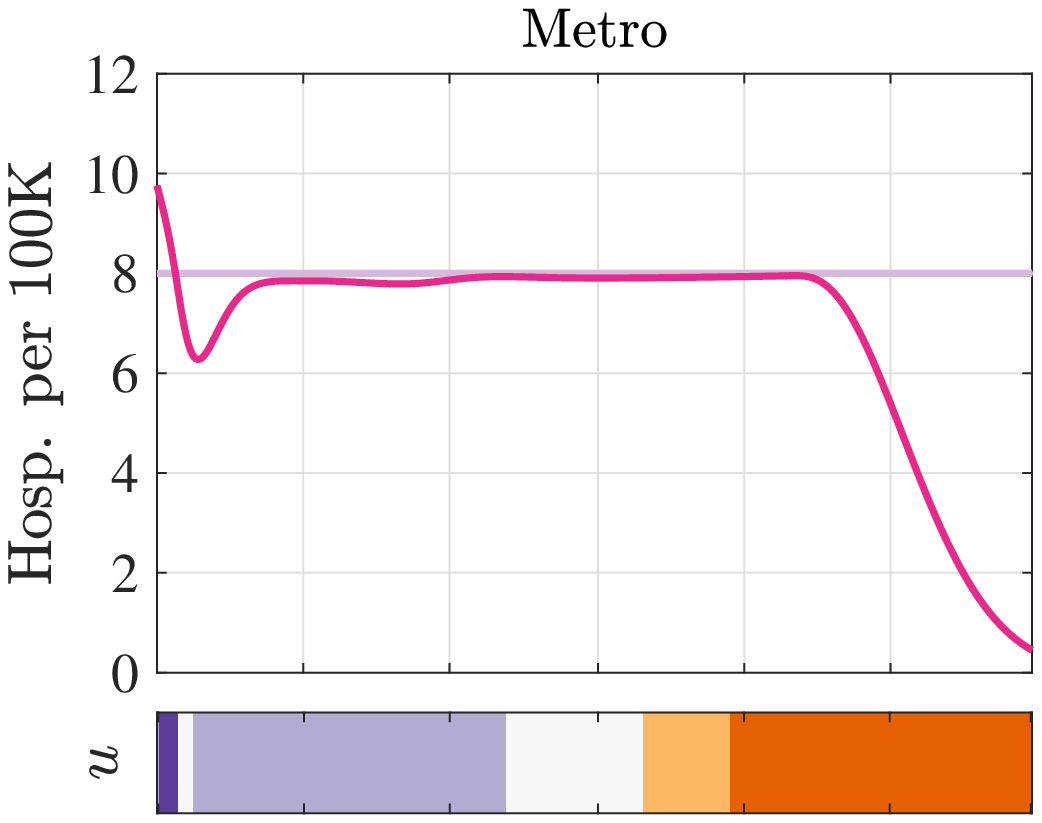}}%
\hspace{.01cm}
\subfigure[]{\includegraphics[width=.24\columnwidth]{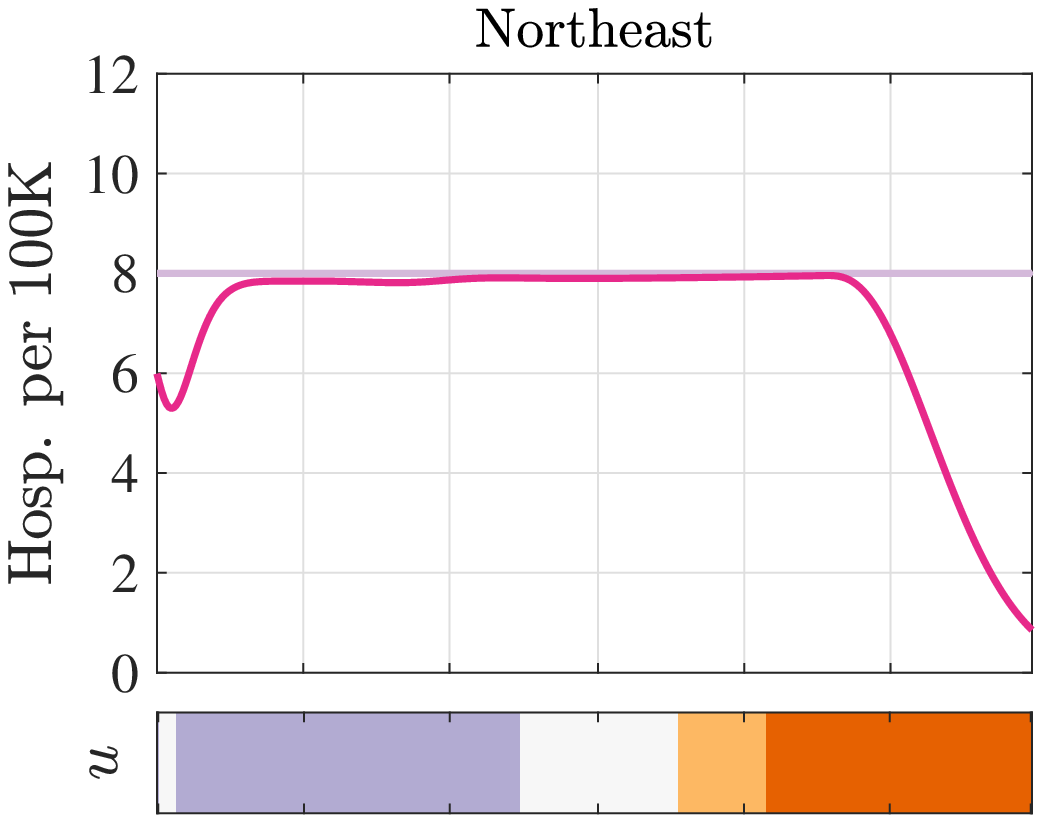}}%
\hspace{.01cm}
\subfigure[]{\includegraphics[width=.24\columnwidth]{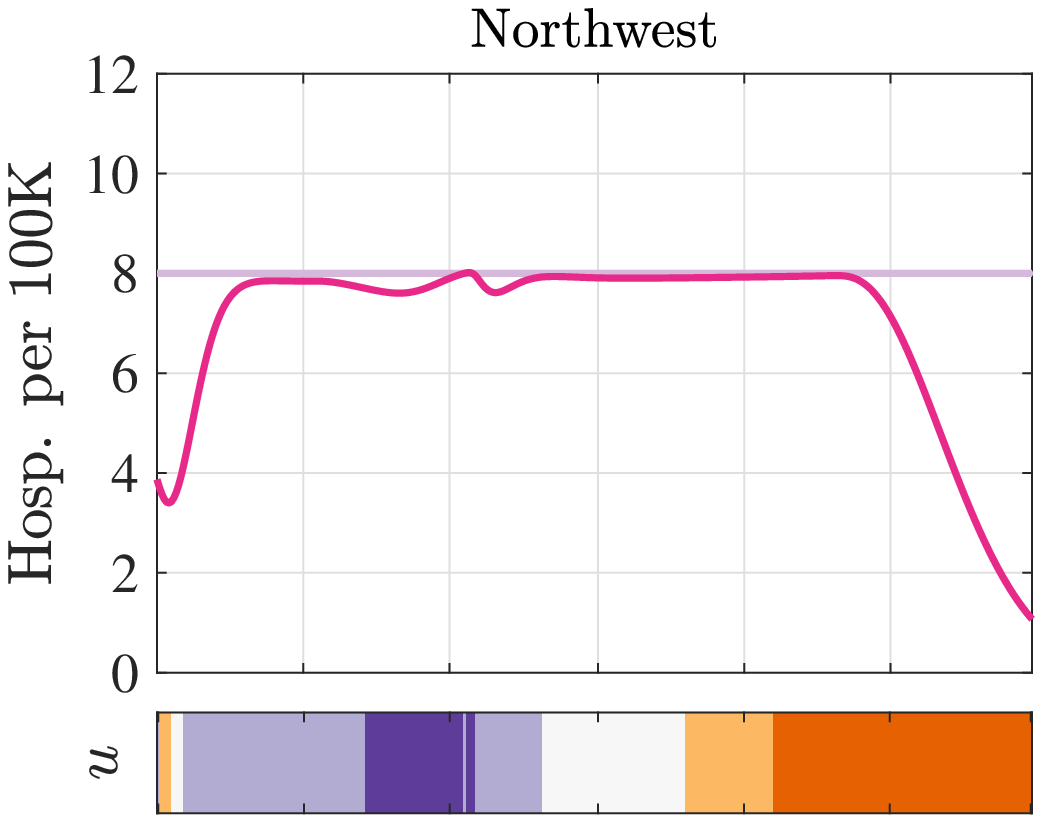}}%
\hspace{.01cm}
\subfigure[]{\includegraphics[width=.24\columnwidth]{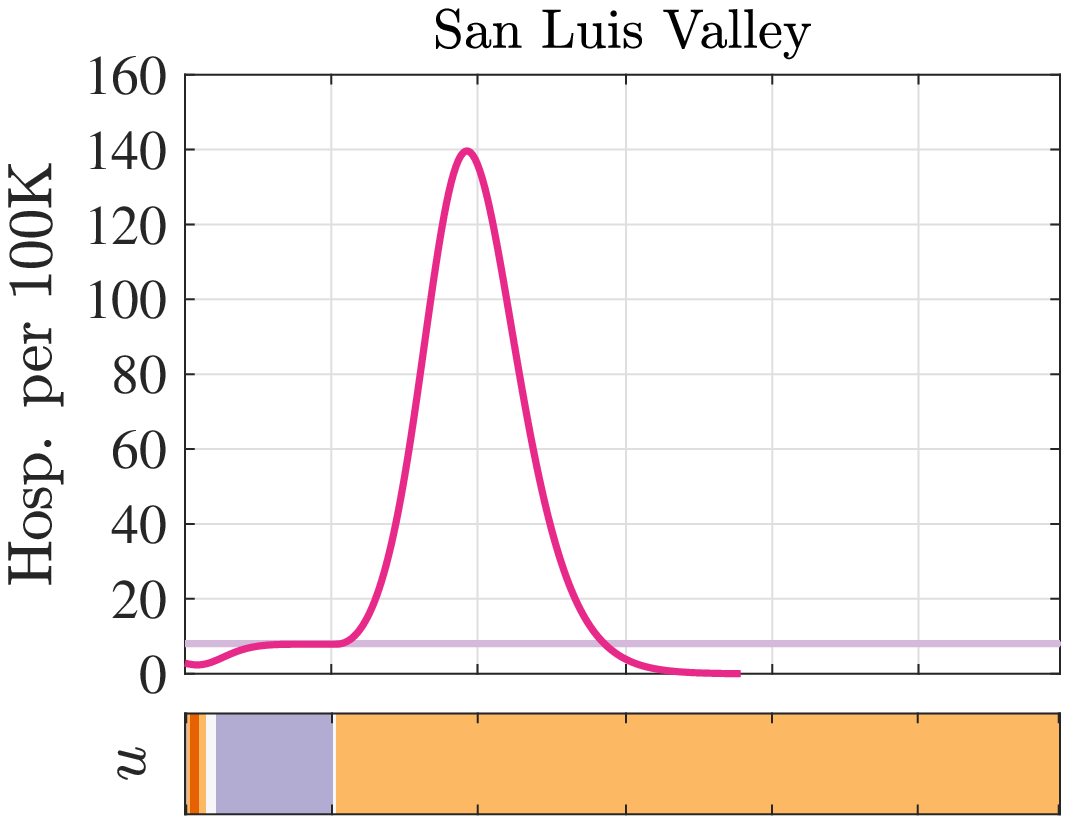}}\\
\subfigure[]{\includegraphics[width=.24\columnwidth]{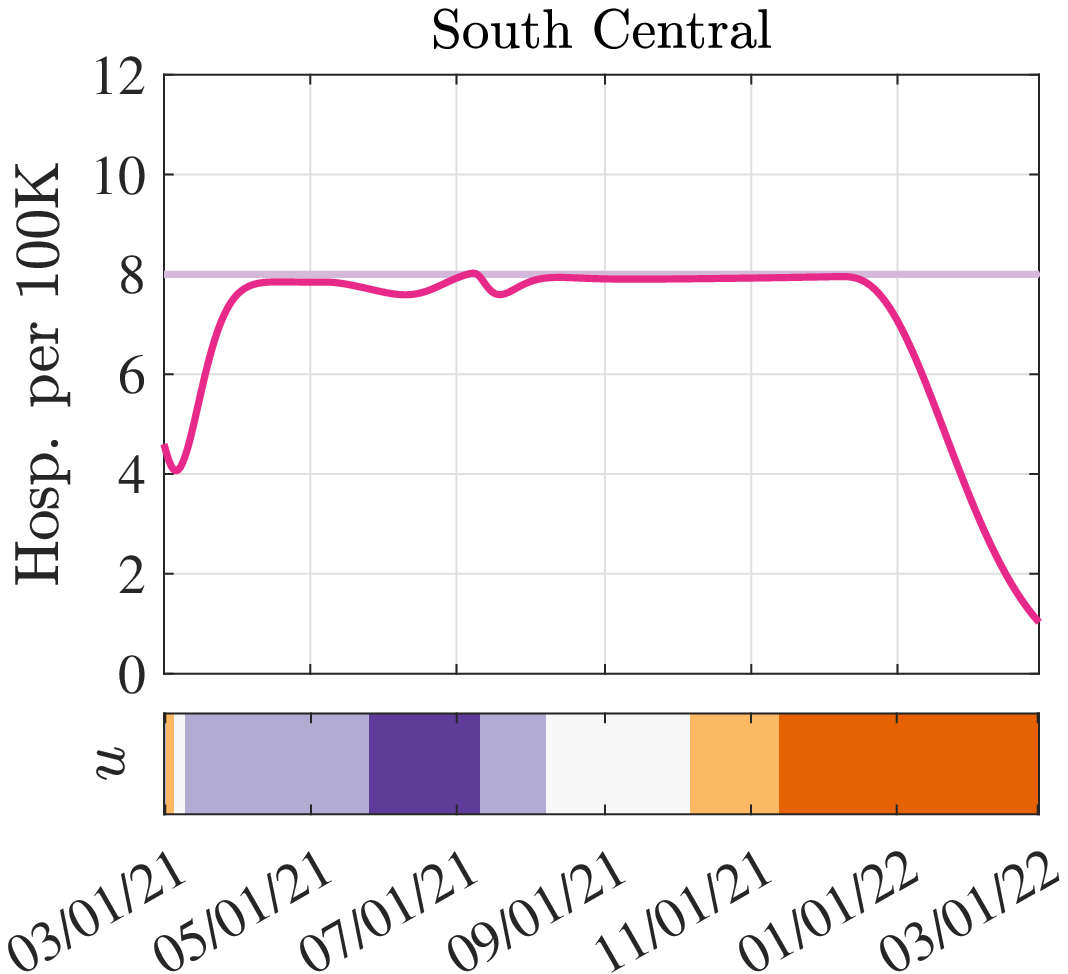}}%
\hspace{.01cm}
\subfigure[]{\includegraphics[width=.24\columnwidth]{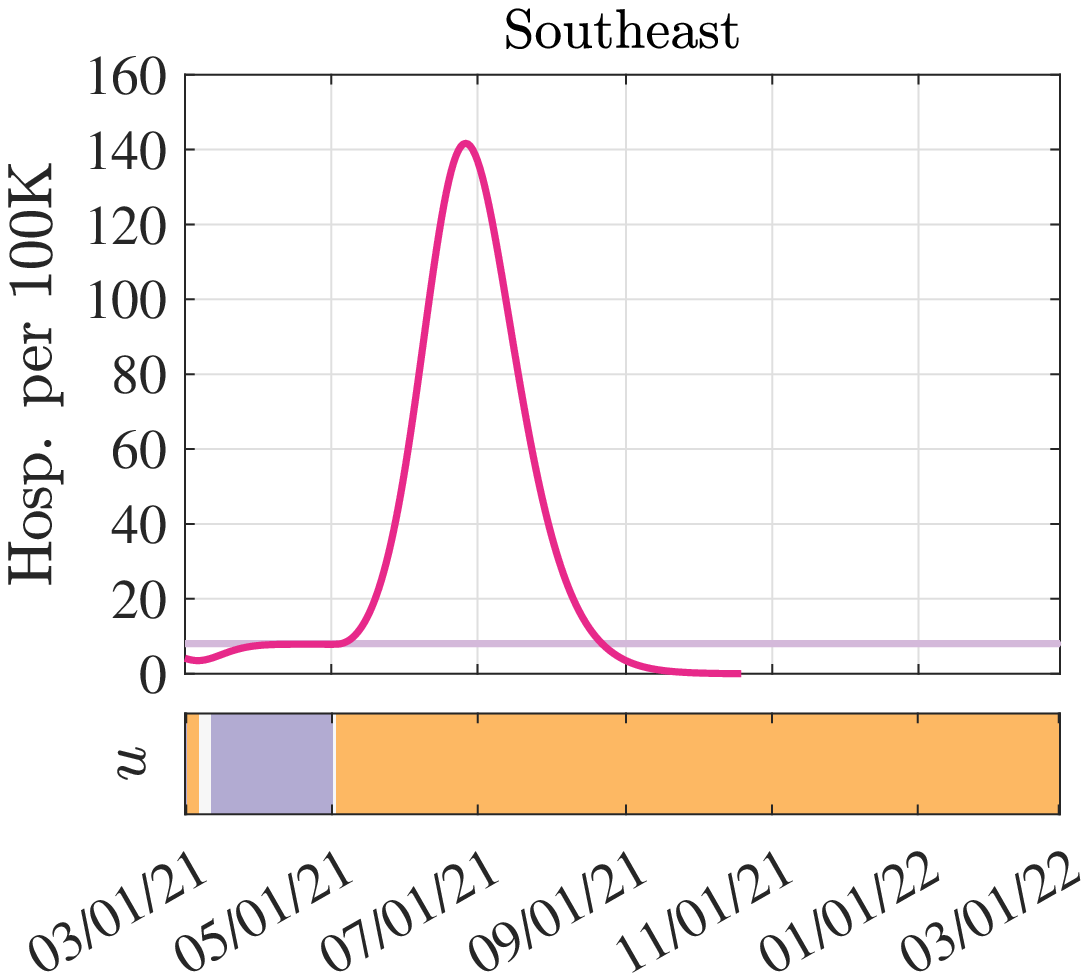}}%
\hspace{.01cm}
\subfigure[]{\includegraphics[width=.24\columnwidth]{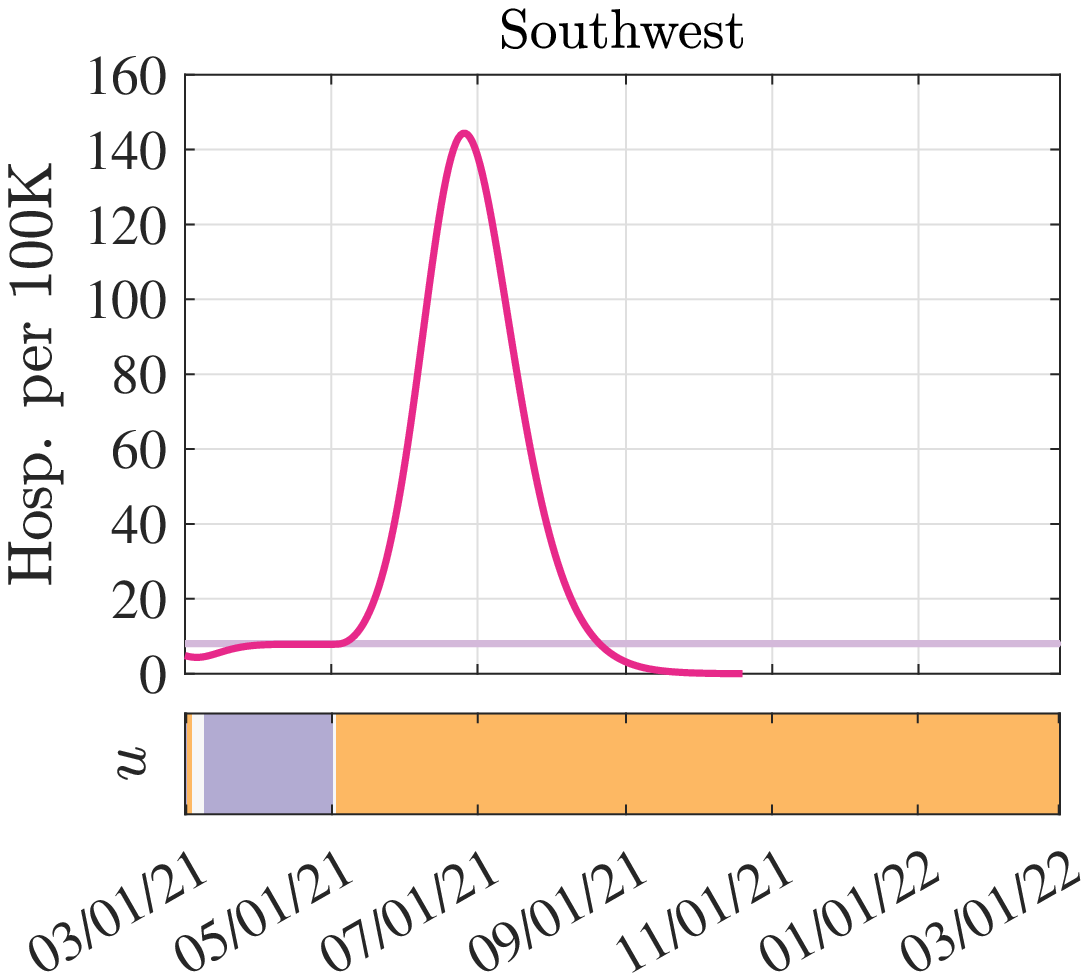}}%
\hspace{.01cm}
\subfigure[]{\includegraphics[width=.24\columnwidth]{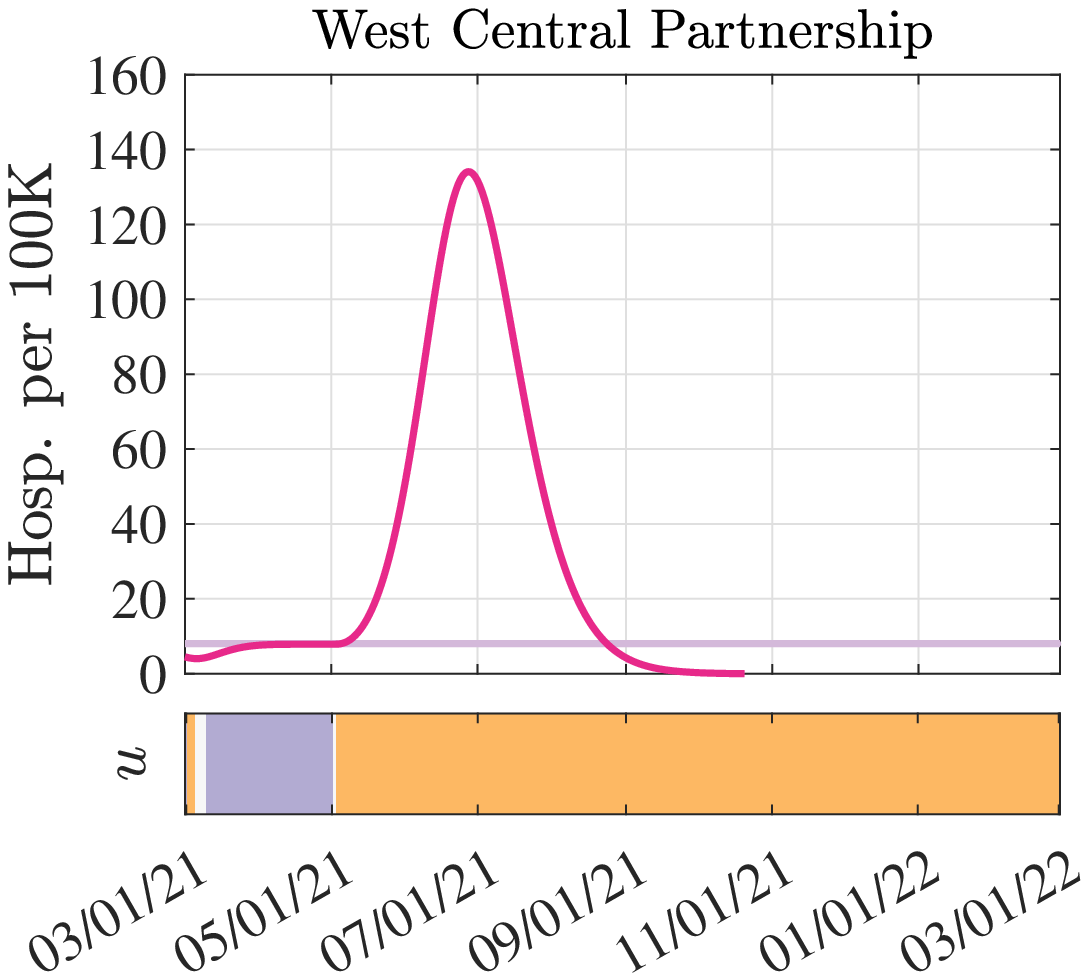}} 
\caption{\emph{\small{
Hospitalizations and controller level over time when regions
with a population of 
$150,000$ people or less (i.e., East Central,
San Luis Valley, Southeast, Southwest, West Central Partnership) drop 
all NPIs on 05/01/21. 
Simulation conducted with vaccination rate $y=20,000$ vax/day.
Each of the 11 panels shows the evolution in time in a different LPHA region 
(see \figurename~\ref{fig:network_regions_graph} for an illustration of the 
connectivity graph).
Solid magenta lines illustrate the evolution of the hospitalized state,
light purple lines show the pre-specified hospitalization limit. 
Heat maps illustrate the required level of NPIs $u_i$, as determined by
the controller. Note that y-scale differs between panels.
}}}
\label{fig:regional2}
\end{figure}

\paragraph{The Value of Coordination.}
\figurename~\ref{fig:regional2} considers a scenario where   
regions with a population of $150,000$ people or less (i.e., East 
Central, San Luis Valley, Southeast, Southwest, and the West Central 
Partnership regions) drop most NPIs as of 05/01/21, returning to $80\%$
contact levels ($u = 0.8$). On this date, the average fraction of fully
vaccinated individuals across the state is $21.29\%$. 
As illustrated by the simulation, such policy will result in a substantial 
violation of the hospitalization limit in all of the five regions that drop 
the NPIs.
Not surprisingly, the regions that decrease NPIs on May 1st, 2021, are 
also the ones that are affected by the highest number of hospitalizations, 
with peaks of over 140 hospitalizations/day per 100K inhabitants around 07/01/21.
Interestingly, our results suggest that in this case three highly-populated 
regions (Central Mountains, Northwest, and South Central) are required to 
decrease the fraction of contact-relevant interactions to below $20\%$
in order to not exceed the hospitalization limit of 8 individuals/day per 
100K inhabitants.
This suggests that high-population regions are at high risk of outbreaks as a result of low levels of control in rural regions. Note that without increasing control in, for instance, the Northwest region, in response to low control levels in surrounding rural areas, the Northwest region would be at risk of large increases in hospitalizations. 
Together, our results indicate that regional heterogeneities can be
exploited by the feedback controller to alleviate the necessary severity
of NPIs to stay below hospitalization limits in interconnected regions, 
however, uncoordinated changes of NPIs in some of the regions will in 
general impact the level of NPIs imposed in all the remaining regions.

We conclude by noting that, although our models are fitted to 
hospitalization data from the state of Colorado, and the simulations in 
figures \ref{fig:regional1}--\ref{fig:regional2} are performed for this
case study, the proposed control framework and feedback control law are
applicable to any geographical region and can be implemented at a 
different geographical granularity.

\section*{Discussion}
\label{sec:discussion}

Vaccinations effectively help to contain the spread of an epidemic  by 
quickly developing immunity in vaccinated individuals, without causing 
severe illness.
With sufficient vaccination rates and high vaccine uptake, NPIs can be 
gradually lifted and social interactions can partially resume in the 
interest of re-establishing economic and social activities. 
Despite optimism over widespread vaccination, a safe return to 
pre-COVID contact behaviors (corresponding to zero NPIs), may still be 
a long way away, dependent on the number of SARS-CoV-2 infections and 
consequent severe COVID-19 cases we are willing to tolerate.
Indeed, if vaccine uptake remains low, policymakers should face the 
possibility of having to either tolerate a high level of NPIs, or a 
high number of severe COVID-19 hospitalizations in the foreseeable 
future.

We examined the conditions under which all NPIs can be safely reverted 
and individuals can resume pre-pandemic contact behaviors.  
We successfully showed that the adopted control method can be used to 
identify necessary increases or decreases in NPIs based on the level of
community risk tolerance or how severe the hospitalization burden can 
be tolerated before governments or individuals will act to impose 
restrictions or change their behaviors.

Control of epidemics is a research area with extensive prior works (see, e.g.,~\cite{RR-RL-CG:09,HB:00,EH-TD:11,DG:88,NG-RR:73,CN-VP-JP:16,LH-DB-TG:04} and pertinent references therein), which include a variety of methodologies that build on model predictive control, model-based  optimal control, periodic lock-downs, etc, or simply heuristics driven by conventional wisdom. Here, we take a novel approach based on online optimization methods for dynamical systems \cite{MC-ED-AB:19,GB-JC-JP-ED:21-tcns,GB-JIP-ED:20-automatica}. Online optimization provides powerful tools to simultaneously control a  dynamical system and steer it to an optimal state configuration, where optimality is quantified according to a pre-specified cost function and constraints that embeds economic and risk tolerance metrics. Model uncertainties constitute the major complexity of the control task  at hand, and call for the development of novel control tools that can  determine how NPIs can be updated over time under limited system  knowledge. The online optimization methods utilized in this work leverages feedback from the  system to adaptively update the control variables in the face of  possible model uncertainties and externalities. Here, the evolution of the pandemic is controlled based on objectives embedded in the optimization problems and by relying on the current number of infectious individuals and a prediction of the function that maps contact levels into number of hospitalizations. The latter can be obtained from data generated by using a transmission model (the approach taken in this paper) or using machine learning tools. The setting investigated here also calls for new theoretical endeavors to uncover the stability properties of networked nonlinear dynamical systems modeling the progression of the epidemics and data-driven controllers, especially when the underlying transmission model accounts for loss of immunity.

In this work, we use a limit on hospitalizations, representing plausible risk tolerance thresholds and stress of the healthcare system, to examine the role of vaccination rate and vaccine uptake on minimum necessary levels of contact. If we assume a low risk-tolerance of no more than $8$ individuals hospitalized/day per 100K inhabitants, 
our results suggest the intriguing possibility that the number of days 
to normality decreases by a factor of two as vaccination uptake 
in the population is increased from $50\%$ to $60\%$.
While increasing vaccination rates will lead to a decreased ``time to 
normal,'' under conservative levels of risk tolerance, safe return to normal may 
not occur until early 2022. Allowing for increased burden of
hospitalizations decreases the time to a safe return to normal behavior, 
but with serious consequences in the form of increased morbidity and 
mortality, even under scenarios with high vaccination rates. 
Vaccine  uptake is a key factor for the return to normal. 
\figurename~\ref{fig:daysTou=1} illustrates the number of days before all NPIs can be lifted in relationship to different levels of hospitalization limits, daily vaccination rate, and maximum vaccine 
uptake. Our results suggest that when the hospitalization limit is maintained below $10$ individuals/day every 100K inhabitants per day, 
at least $300$ days are necessary to lift all NPIs (given  a 
vaccination rate of 15,000 vaccines/day). When the vaccination uptake is as low as $40\%$, this prediction increases to about  $700$ days.  
\figurename~\ref{fig:daysTou=1}(c) illustrates the cumulative number of 
deaths from 03/01/21 for different levels of allowed hospitalizations
and vaccination rates. Our results suggests that the number of deaths 
grows linearly with the number of allowed hospitalizations, independent of the rate or uptake of vaccinations.
This behavior is due to model assumptions, such that the number of deaths is 
proportional to the number of infections and of hospitalized
individuals, which are maintained constant over time  by the  control  method. In practice, this assumption does not hold as vaccinations decrease the death rate among hospitalized individuals proportional to the high-risk population vaccinated, but this model does not account for this heterogeneity.  

Regional heterogeneity complicates this picture. Even when only  low-population regions with relatively low contact rates may begin 
to return to relative normality far more quickly, mobility across regions plays a key role. Due to people moving and interacting across regions, removing NPIs too quickly even in regions of low population 
density can still lead to dire consequences in nearby high-density 
regions. In our interconnected world these findings can be generalized to both smaller and greater spatial scales; the application of our method can unveil important intrinsic  dependencies that should be fully taken into account to effectively control the spread of the infection.

We acknowledge that our findings come with a number of relevant 
limitations. First, the outcomes are dependent on numerous assumptions 
about baseline transmission rate, probability of hospitalization, and 
other parameter values estimated from previous modeling studies that may 
impact our results. Second, we chose not to account for age or the 
differences between asymptomatic and symptomatic transmission for 
simplicity. Third, despite accounting for regional heterogeneity in contact 
rates and baseline transmission, superspreader events and smaller 
non-homogeneous spatial units play a large role at this stage in the 
pandemic. Compounding this, vaccine distribution is occurring in a manner which reinforces 
pre-existing health disparities, due to issues of both access and 
hesitancy~\cite{ndugga2021latest}. This creates pockets of high-risk 
unvaccinated populations, which are sufficient to sustain transmission, 
even with high vaccination rates overall. Our model cannot account for 
this type of clustering of behavior or risk, which are important in 
understanding the probability of achieving sufficiently low levels of 
SARS-CoV-2 transmission. 
Fourth, we do not account for the recent introduction and proliferation 
of numerous variant strains which have the potential to substantially 
alter transmission dynamics and vaccine efficacy~\cite{fontanet2021sars}; future work will account for variants and include them in the proposed methodology. Over time, current vaccines may be less effective at preventing infection
due to new circulating variants, preventing attainment of herd immunity 
even with high rates of vaccination uptake. Similarly, the duration of 
immunity obtained from current vaccines against currently circulating 
variants is unknown. For the purpose of this study, we assumed an 
effectiveness of two years (which may be overly optimistic). If the 
duration of vaccine-derived immunity is shorter in practice, given 
feasible vaccination rates, complete relaxation of NPIs may never be 
attainable~\cite{aruffo2021community}. When better data on  
duration of immunity become available, the model utilized here can be modified
accordingly. Finally, we also acknowledge that it may be difficult to 
directly translate a variable between $0$ and $1$ into precise actions 
such as mask mandates, school closures, and business capacity limits, 
and especially personal decisions; however, current studies are looking 
at regressing the reproduction number of COVID-19 against different 
NPIs~\cite{liu2021impact}, and this can naturally be connected to the 
variable $u$ of the controller. 

Our findings are in agreement with previous modeling studies which have 
stressed the need to maintain current levels of NPIs and decreased 
contact for the near future, even in the context of current vaccination 
strategies~\cite{moore2021vaccination,li2021returning,love2021continued}.
Several recent studies have also questioned whether herd immunity through
vaccination is achievable at all, given the current vaccines available 
and the high prevalence of vaccine hesitancy~\cite{zachreson2021will}.
Given these factors, the possibility has been raised that SARS-CoV-2 will
become an endemic virus circulating regularly in the 
population~\cite{lavine2021immunological} and our concept of a ``return to
normal'' will have to be reframed. Currently (04/09/21), in the state of
Colorado, transmission has flattened with approximately 400 individuals 
currently hospitalized. Vaccination rates are dropping rapidly and whether or 
not we can reach $70\%$ vaccine uptake is uncertain.

\section*{Methods}
\label{sec:methods}

The evolution of the epidemic is modeled by using a 
Susceptible-Exposed-Infectious-Hospitalized-Recovered-Vaccinated-Susceptible (SEIHRVS) compartmental model. We begin by illustrating the 
single-region model, we then extend the model to account for regional
heterogeneities, and lastly we illustrate the control method.

\begin{figure}[t]
\centering
\includegraphics[width=10cm]{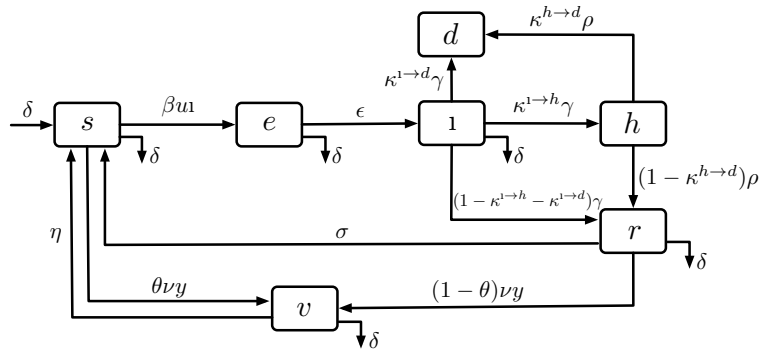} 
\caption{\emph{\small Block diagram of the  compartmental model adopted to generate data. In particular, the model is for a single-region. Model equations and extensions to the multi-region model are discussed in Methods.}}
\label{fig:compartmentalModel}
\end{figure}

\paragraph{Single-Region Modeling.} For the single-population case, we 
utilize a  variation of the  Susceptible-Exposed-Infectious-Recovered 
(SEIR) model \cite{JA-PV:03} that
accounts for hospitalizations, vaccinations, and loss of immunity. In  
particular, we consider a transmission model with states:
Susceptible ($s$), Exposed ($e$) Infectious ($\i$), Hospitalized ($h$),  
Recovered ($r$), Vaccinated ($v$), and Deceased ($d$). Infectious
individuals can infect susceptible ones with a transmission rate 
$\beta > 0$. To model NPIs, we let $u \in [0,1]$ be a scalar parameter 
that  specifies the level of permitted social activity (or contact 
levels within the region). 
The special case $u = 0$ models a full lock-down, while $u = 1$ 
corresponds to ``zero'' NPIs (and hence, a return to pre-pandemic contact 
levels). Noticing that different levels of NPIs result in different 
transmission-relevant  contact levels, the overall model of epidemic 
transmission is given by the following differential equations:
\begin{align}
\label{eq:SEIR}
\dot s &= - \beta u s \i - \theta \nu y -  \delta s + \delta + \sigma r + \eta 
v,\nonumber\\
\dot e &=  - \epsilon e - \delta e + \beta u s \i,\nonumber\\
\dot \i &= - \gamma \i - \delta \i + \epsilon e,
\nonumber\\
\dot h &=  - \rho h  + \kappa^{\i \rightarrow h} \gamma \i,\nonumber\\
\dot r &= -\sigma r - \delta r - (1-\theta ) \nu y
+ (1- \kappa^{\i \rightarrow h} - \kappa^{\i \rightarrow d})\gamma \i
+ (1- \kappa^{h \rightarrow d})\rho h,\nonumber\\
\dot v &= - \eta v - \delta v + \nu y \nonumber,\\
\dot d &= \kappa^{\i \rightarrow d} \gamma \i + \kappa^{h \rightarrow d} \rho h,
\end{align}
where we use $\dot{x}:=\frac{d}{dt}x(t)$ to denote the time-derivative
of a scalar-valued variable $x(t)$ that is a function of time $t$.
In \eqref{eq:SEIR}, whenever $u \in [0,1)$, the effective transmission rate is
reduced from $\beta$ to $\beta u$  according to the imposed NPIs. 
Individuals become infectious after an incubation period $1/\epsilon$, and they recover at a rate $\gamma > 0$. After being infectious, a fraction of individuals 
$\kappa^{\i \rightarrow d} \in[0,1]$ dies, and a fraction 
$\kappa^{\i \rightarrow h} \in[0,1]$ is hospitalized.
The fraction $1-\kappa^{\i \rightarrow d}- \kappa^{\i \rightarrow h}$ 
quantifies the individuals who recover without hospitalization.
Hospitalized individuals recover at a rate $\rho>0$.
After being hospitalized, a fraction $\kappa^{h \rightarrow d} \in[0,1]$
of individuals die,  while $1-\kappa^{h \rightarrow d}$ recover.
Recovered individuals lose immunity at a rate $\sigma>0$, thus 
returning in the susceptible compartment. We denote by $y$ the 
vaccination rate and by  $\nu$ the vaccination efficacy.
Individuals are vaccinated regardless of their prior infection history, 
and we let $\theta\in[0,1]$ be the fraction of vaccines that is 
administered to individuals in the $s$ compartment, while  $(1-\theta)$ 
vaccines are administered to individuals in the $r$ compartment. 
Finally, $\delta>0$ describes the population birth/death rate. The 
compartmental model corresponding to the differential 
equations~\eqref{eq:SEIR} is illustrated in 
\figurename~\ref{fig:compartmentalModel}. 

\paragraph{Multi-Region and Mobility Modeling.}
We consider a model of disease transmission that is organized into a 
group of geographical subregions, where individuals make short-term 
(e.g.~daily) inter-regional movements or transits. The assumption 
that transits are short-term models scenarios where exposed 
individuals return to the corresponding region of residence before, 
eventually, becoming infectious. 
\figurename~\ref{fig:network_regions_graph} illustrates the 
partitioning of the state of Colorado according to Local Public 
Health Agency (LPHA) regions, and illustrates the flow of mobility 
between regions.

To model contact-relevant interactions among residents of different 
regions, we adopt a graph $\mc G = (\mc V, \mc E)$ where 
$\mc V = \until N$ denotes the set of nodes (regions), and 
$\mc E \subseteq \mc V \times \mc V$ denotes the set of edges 
(links between regions).  We model the coupling between regions by 
assuming that a fraction $a_{ij} \in [0,1]$ of residents of region $j$ 
travel to region $i$ and interact with its residents, 
with $\sum_{j} a_{ij}=1$
(here, $a_{ij}$ 
denotes the intensity of contacts between residents of region $i$ 
and residents of region $j$ at region $i$). To model the epidemic spread in 
the regions, we assume that each individual that is a resident of region 
$i$ is categorized into one of the seven compartments $s_i$, $e_i$, 
$\i_i$, $h_i$, $r_i$, $v_i$, $d_i$ (where the subscript $i$ denotes the 
index of the region). Similarly to \eqref{eq:SEIR}, each state of 
\eqref{eq:SEIR_network} represents the fraction of individuals in the 
corresponding compartment, with $s_i+e_i+\i_i+r_i+h_i+v_i+d_i=1$ at all 
times. The extension of~\eqref{eq:SEIR} to a multi-region setting is given by:
\begin{align}
\label{eq:SEIR_network}
\dot s_i &= - \sum_{j=1}^N \beta u_i a_{ij} s_i \i_j - \theta_i \nu_i y_i
-  \delta_i s_i + \delta_i + \sigma_i r_i + \eta_i v_i,\nonumber\\
\dot e_i &=  - \epsilon_i e_i - \delta_i e_i 
+ \sum_{j=1}^N \beta u_i a_{ij} s_i \i_j,\nonumber\\
\dot \i_i &= - \gamma_i \i_i - \delta_i \i_i + \epsilon_i e_i,
\nonumber\\
\dot h_i &=  - \rho_i h_i  + \kappa^{\i \rightarrow h}_i \gamma_i \i_i,\nonumber\\
\dot r_i &= -\sigma_i r_i - \delta_i r_i - (1-\theta_i ) \nu_i y_i
+ (1- \kappa^{\i \rightarrow h}_i - \kappa^{\i \rightarrow d}_i)
 \gamma_i \i_i
+ (1- \kappa^{h \rightarrow d}_i)\rho_i h_i,\nonumber\\
\dot v_i &= - \eta_i v_i - \delta_i v_i + \nu_i y_i \nonumber,\\
\dot d_i &= \kappa^{\i \rightarrow d}_i \gamma_i \i_i 
+ \kappa^{h \rightarrow d}_i \rho_i h_i.
\end{align}
We note that 
\eqref{eq:SEIR_network} allows for different levels of NPIs across the 
regions, where the variable  $u_i \in [0,1]$ describes the permitted  level 
of transmission-relevant contacts within region $i$. 
This model is motivated by recent works that found the 
regional progress of influenza much more correlated with the movement of
individuals rather than geographic distances 
\cite{CV-OB-DS-LS-MM-BG:06}. For these reasons, network models are 
widely used in the literature to take into account spatial propagation 
effects \cite{LS-KD:95,MT-PB-AD:14,JA-PV:03,AT-TP-NG-etal:20,QM-YYL-AO:20,mesbahi2010graph}. \\

\noindent\textbf{Defining a Metric for Control}. For the purpose of this study we utilize a single control variable per region to encompass all policy measures and behavioral changes in a region that decrease transmission of SARS-CoV-2. This variable captures the impact of several policy measures: including mask mandates, school closures, business capacity limits, as well as personal decisions such as hand-washing, mask wearing, and moving socialization outside that otherwise would have occurred inside. 
At present, it is practically intractable to disentangle the specific impact that individual interventions have, as multiple complex interventions are introduced simultaneously and the population is reacting in a continuous manner to changing risk perception influenced by divergent policy and messaging at the local, state, and federal level~\cite{SF-SM-AG-etal:20}. While NPIs implemented by policy-makers are important, a large portion of transmission reducing behaviors are a result of individual-level risk assessment and behaviors, in response to perceived community transmission~\cite{yan2020measuring}. Thus, even as policy makers begin to relax NPIs at the state and regional levels, individuals will continue to make decisions based on their own perceptions of risk, which are directly impacted by hospitalization and infection levels in the community. As a result of these effects, we recognize that an actual implementation of various NPIs may have a high variance. For example, even if all policy measures are lifted, as prevalence remains high, it might be unlikely that some  individuals will return to pre-pandemic contact behavior. Likewise, even at the height of restrictions, when stay at home orders were in place, it may not be  possible to control transmission entirely given the necessity of ongoing essential work, grocery shopping, etc (therefore, $u$ can approach $0$, but cannot be set to $u = 0$ in practice), or because of possible violations of restrictions.

\subsection*{Feedback optimization theory for NPI}

For the design of feedback controllers, we take an approach inspired by 
recent advances in  feedback optimization of dynamical systems~\cite{brunner2012feedback,MC-ED-AB:19,GB-JC-JP-ED:21-tcns,hauswirth2020timescale,GB-JIP-ED:20-automatica} and, in particular, we develop a new data-driven optimization approach based on the analytical framework 
described in~\cite{GB-JC-JP-ED:21-tcns}. 

\paragraph{Formalizing risk tolerance and social objectives.} The proposed technical approach builds upon formulating an optimization problem that captures the  desired social and economic metrics, and incorporates constraints  related to risk tolerance. 

To this end, let $\sbs{h}{lim,$i$}$ denote the maximum allowable 
number of daily hospitalized individuals in region $i$, and let $\map{\phi_i}{[0,1]}{\realpos}$ be a function of the  decision variable $u_i$ that models the societal impact induced by  the introduction of NPIs in region $i$~\cite{QM-YYL-AO:20}; this includes the economic impact of the raised control measures, and/or the societal response to restriction orders. 
Similarly, let $\map{\psi_i}{[0,1]}{\realpos}$ be a function of the number of 
infectious individuals $\i_i$ that models the societal losses due to a high number of
infections in region $i$; this includes e.g. the cost of hospitalizations.
Mathematically, we assume that $u\mapsto\phi_i(u)$ is a differentiable 
function~\cite{nocedal2006numerical} (see Supplementary Information). Since societal losses, in general, do not increase as NPIs are lifted, we also assume that $u\mapsto\phi_i(u)$ is a non-increasing function in its domain. For additional remarks on the cost of NPIs we refer the reader to, e.g., \cite{QM-YYL-AO:20,chen2020tracking,ugarov2020inclusive}. To capture the relationship between NPI variables $u = (u_1, u_2 \ldots, n_N)$ and the number of infectious individuals at the endemic equilibrium (i.e. when $t = \infty$), we denote by 
$u \mapsto \mc F_i(u)$ the function that maps the instantaneous values of NPIs $u$ to the fraction of infections in the $i$-th region at the endemic equilibrium. Hospitalizations can be naturally related to the instantaneous level of NPIs as well as to the current state of the pandemic. To this end, we let $x_i:=(s_i, e_i, \i_i, h_i, r_i, v_i, d_i)$ 
denote the joint epidemic state in region $i$, and 
$x:=(x_1, x_2, \dots , x_N)$ denote the joint epidemic state of the 
network.  We denote by $u \mapsto \mc H_i(u; x)$ the function that maps 
the instantaneous level of NPIs $u$ into the peak of hospitalizations in 
the region $i$, given the current state $x$ of the model. 

With these definitions in place, we formulate the problem of 
optimizing the choice of NPIs while guaranteeing that the number of 
hospitalizations remains below the pre-specified limit at all 
times as follows:
\begin{align}
\label{opt:NPIandHospitalizationFormulationNet}
\min_{\{u_1, \ldots, u_N\}}~~~ & \sum_{i = 1}^N \phi_i(u_i) + \psi_i(\mc F_i(u)) & & \text{(cost depends 
on NPI policy and infections in each area)} \nonumber\\
\text{s.t.}~~~  &  \mc H_i(u; x) \leq \sbs{h}{lim,$i$}, \,\,\,\,\, i = 1, \ldots, N, & &\text{(maximum allowable hospitalizations in area $i$)}\nonumber \\
& u_i \in [0,1], \hspace{1.4cm} i = 1, \ldots, N. & & \text{(feasible NPI policy for area $i$)}
\end{align}
Solutions $u^*$ of the above optimization problem describe a level of NPIs
that balances economic and social costs, while ensuring that
hospitalizations limits are not exceeded in each region. In the following,
we develop an optimization-based feedback control method to continuously 
calibrate the permitted level of social interactions $u$ to meet the 
criteria set forth in the optimization 
problem~\eqref{opt:NPIandHospitalizationFormulationNet}. We also  
remark that~\eqref{opt:NPIandHospitalizationFormulationNet}
includes a single-region model as a subcase (see Supplementary Information).

\begin{figure}[t]
\centering \includegraphics[width=16.0cm]{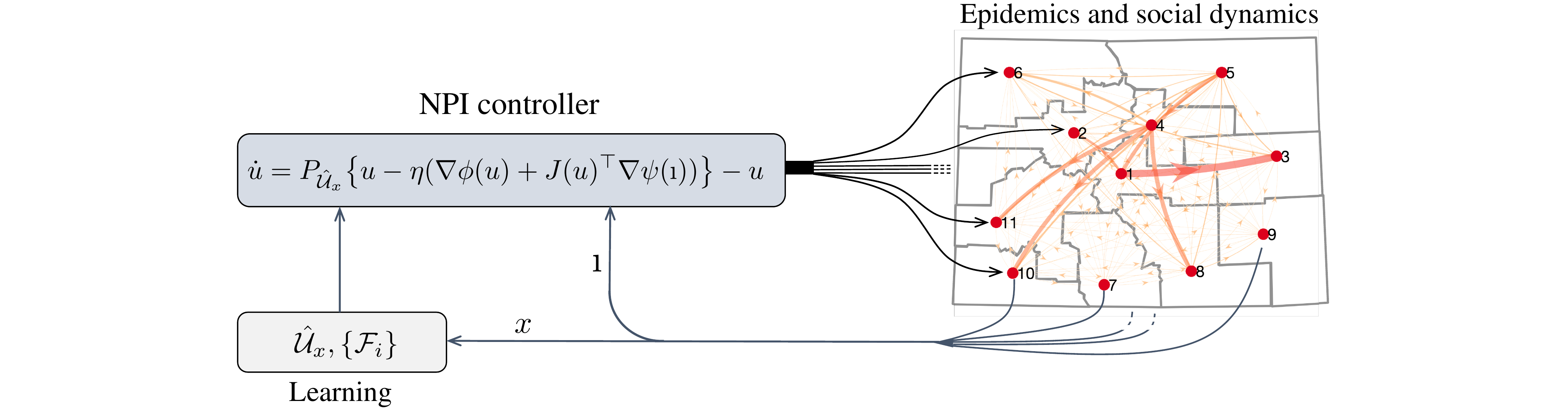}
\caption{\emph{\small{Implementation of the NPI controller. The example refers to the state of Colorado, where each region represents a Local Public Health Agency. }}}
\label{fig:controller_co}
\end{figure}

\paragraph{Feedback controller design for NPIs.} We begin by defining the 
set of feasible NPIs, as described in the  optimization problem 
\eqref{opt:NPIandHospitalizationFormulationNet}. 
For a fixed epidemic state $x$, the feasible region of 
\eqref{opt:NPIandHospitalizationFormulationNet} is given by:
\begin{align*}
\mc U_x := \setdef{ u = (u_1, \dots u_n) }{
u_i \in [0,1], 
\mc H_i(u; x) - \sbs{h}{lim,$i$} \leq 0, 
\text{ for all } i = 1, \dots N} \, .
\end{align*}
We note that, because the input-to-peak of hospitalization map 
$\mc H(u; x)$ is parametrized by the instantaneous state of the system $x$, the set $\mc U_x$ is also parametrized by $x$.   
When the set $\mc U_x$ is non-convex, we consider a convex approximation $\hat{\mc U}_x$ as explained in the Data-driven implementation section. 
With this definition, a function $t\mapsto u(t)$ for the NPIs can be 
obtained as a solution of the following dynamical system:
\begin{align}
\label{eq:controllerGradientFlow}
\dot u &= P_{\hat{\mc U}_x} \big(u - \eta (\nabla \phi(u) 
+ J(u)^\top \nabla \psi(\i))\big) - u,
\end{align}
where $\phi(u) := \sum_{i = 1}^N \phi_i(u_i)$ and $\psi(i) := \sum_{i = 1}^N  \psi_i(i_i)$ for brevity, $J(u)$ is the Jacobian matrix collecting $\{\partial_u \mc F_i (u)\}$,  $\eta>0$ is a tunable parameter of the controller, and $P_{\mc U_x}$ denotes the Euclidean projection operator; namely, given $z \in \real^n$ and a convex set $\mc U \subseteq \real^n$,
\begin{align*}
P_{\mc U} (z) = \arg \min_{w \in \mc U} \norm{w-z}.
\end{align*}
We note that the optimization-based controller  \eqref{eq:controllerGradientFlow} 
leverages two types of feedback: (i) it uses the instantaneous fraction of infectious 
individuals $\i$, and (ii) it relies on a projection onto the set $\hat{\mc U}_x$, which 
is parametrized by the instantaneous state of the system.
For these reasons, the control dynamics \eqref{eq:controllerGradientFlow} describe a \emph{dynamic state-feedback controller} for the NPIs. Critically, the controller relies on the knowledge of the maps $u \mapsto \mc F_i(u)$ and $u \mapsto \mc H_i(u; x)$. These maps are estimated from data, as explained in ``Data-driven implementation''. An illustrative example of the implementation of the  controllers~\eqref{eq:controllerGradientFlow} is provided in \figurename~\ref{fig:controller_co}.

\paragraph{Local implementation.}  Due to the coupling introduced by the dependence of functions $\mc H_i$ 
and $ \mc F_i$ on the (entire) vector of control variables $u$, the 
implementation of the optimization-based feedback controller \eqref{eq:controllerGradientFlow} critically requires full knowledge of the state, control vector $u$, and of the (gradients of) the cost functions  $\phi_1, \dots \phi_N, \psi_1, \dots ,\psi_N$. 
Therefore, it requires an implementation in a central location (for example, at the state level in the example in \figurename~\ref{fig:network_regions_graph}). When this implementation is not feasible, we consider an approximation of the functions
$\mc H_i$ that accounts only for the effects of the local NPI policies  in area $i$, namely, we approximate the value of the peak of 
hospitalizations $\mc H_i(u; x)$ by 
$\tilde{\mc H}_i(u_i; u_{-i}, x)$, where 
$u_{-i} = \{u_1, \ldots, u_{i-1}, u_{i+1}, \ldots, u_N \}$ is 
treated as a constant parameter.
By using this approximation, we redefine the set of feasible NPI 
policies in subregion $i$  as:
\begin{align*}
\hat{\mc U}_{x,i} := \setdef{ u}{u_i \in [0,1] \text{ and } 
\tilde{\mc H}_i(u_i; u_{-i}, x) - \sbs{h}{lim,$i$} \leq 0},
\end{align*}
where the 
map $\tilde{\mc H}_i$ is obtained numerically.
By using this approximation, the  distributed controller reads 
as:
\begin{align}
\dot u_1 &= P_{\hat{\mc{U}}_{x,1}}[u_1 - \eta (\partial \phi_1(u_1)
+ \partial_{u_1} \mc F_1(u) \partial \psi_1(\i_1))] - u_1, \nonumber \\
\dot u_2 &= P_{\hat{\mc{U}}_{x,2}}[u_2 - \eta (\partial \phi_2(u_2)
+ \partial_{u_2} \mc F_2(u) \partial \psi_2(\i_2))] - u_2, \nonumber \\
\vdots  & \hspace{4.0cm} \vdots \nonumber \\
\dot u_N &= P_{\hat{\mc{U}}_{x,N}}[u_N - \eta (\partial \phi_N(u_N)
+ \partial_{u_N} \mc F_N(u) \partial \psi_N(\i_N))] - u_N\, .
\label{eq:controllerGradientFlow_approximate2}
\end{align}
Each region $i$ can update its NPI policy $u_i$ locally by only relying on the  knowledge of: (i) the current fraction of infectious individual $\i_i$ in  
region $i$, (ii) the current value of NPI in the network $u_{-i}$, (iii) 
an estimate of the partial derivative $\partial_{u_i} \mc F_i(u)$, which 
can be estimated locally at each region by simulating the the 
model~\eqref{eq:SEIR_network}, and (iv) the feasible set $\mc{U}_{x,i}$. 
We note that, although the controllers are implemented locally in each region, coordination between regions naturally emerges because of the connectivity in the SEIHRVS model. 

\paragraph{Data-driven implementation.}  

The maps $\mc H_i$ and $\mc F_i$  can be estimated from data using 
function estimation methods. In particular, to estimate these maps 
we simulated the dynamics \eqref{eq:SEIR} with initial conditions 
set equal to the instantaneous state of the model, and for different  
values of the 
control parameters, chosen in a neighborhood of the current value of
$u$. For a set of fixed values for the control parameters, we 
simulated the dynamics \eqref{eq:SEIR} and obtained the values of 
the peak hospitalizations and the infections at steady state. With 
these values, we utilized function estimation methods to obtain the 
maps $\mc H_i$ and $\mc F_i$ (see Supplementary Information). We 
note that that the estimated map $\hat{\mc H}_i$ is required to be 
convex in order to build a feasibility set $\hat{\mc U}_x$ that is 
convex (which is important in order to have a well-defined 
projection in our controller)~\cite{boyd2004convex,nocedal2006numerical}.

\subsection*{Model Fitting and Data Acquisition}

\paragraph{Model fitting from data.}
We organized the model-fitting phase into two main stages. First, we
fitted the SEIHRVS model by combining model parameters from 
\cite{AB-EC-DG-etal-b:21} with hospitalization data, and we used a 
prediction-correction algorithm to minimize the fitting error.
The fitted model parameters used in our simulations are reported in 
Table \ref{tab:modelParameters}. 
Second, we used cell-phone usage from SafeGraph to estimate the 
interaction matrix of the network. 
The travel volume from an origin region to destination region on a given date is calculated using the Destination Census Block Groups (CBGs) metric in the social distancing data provided by the Safegraph COVID-19 Data Consortium. The destination CBGs metric is defined as: The number of devices with a home in [a CBG in origin region] that stopped in [a CBG in destination region] for $>1$ minute [on a given day]. The ``home'' of the device refers to the most common nighttime location for the device over the prior six weeks. 
The share activity in region $j$ coming from region $i$, used in this analysis, is the activity in $j$ from $i$ divided by the total activity in $j$ across all origins -- where ``activity in $j$ from $i$'' is the sum of the destination CBGs metric for all origin CBGs in $i$ and all destination CBGs in $j$". The summation across CBGs does not perform any deduplication, and so the total activity does not represent unique devices on a given day. Instead, it can be interpreted as the total visits to CBGs in the destination region by devices from the origin region, counting a device that stopped for at least a minute in two different CBGs as two visits. This approach has advantages and disadvantages. It may weight visitors who stay longer and move around during their visit -- going to restaurants, parks, and other locales that are not within the same CBG as their hotels -- more heavily than visitors whose stay is brief or who limit their movement to a small area. The advantage is that longer stays and more movement tend to carry more risk of COVID-19 transmission, and so it helps us capture the impact of restrictions on travel-induced COVID risk. 
Additional information on the data is provided in the Supplementary Information.

\paragraph{Data sources.}

\emph{Mobility data} was obtained from publicly-available Safegraph 
data (\url{https://docs.safegraph.com/docs/social-distancing-metrics}).

\emph{Hospitalization data} was obtained from the Colorado Department 
of Public Health and Environment and came from two different 
data sources. Colorado state-wide daily hospitalization census 
was used to fit the single-region model and came from 
EMResources. This data is available publicly online at 
\url{https://covid19.colorado.gov/data}. Regional-level daily 
hospitalization census data was used to fit the multi-region model and was obtained from COVID Patient Hospitalization Surveillance, data is posted publicly at \url{https://github.com/agb85/covid-19/tree/master/Regional%20Models}.

 \emph{Parameter Values} were obtained from previous modeling
 works \cite{AB-EC-DG-etal-b:21,AB-EC-DG-etal:21} or fit to data and are outlined in detail in the Supplementary Information.

\section{Data and Code Availability}
Mobility data was obtained from Safegraph, publicly available at the 
website: 
\url{https://docs.safegraph.com/docs/social-distancing-metrics}.
Hospitalization data was obtained from the Colorado Department 
of Public Health and Environment. Colorado state-wide daily 
hospitalization was obtained from EMResources, publicly 
available at \url{https://covid19.colorado.gov/data}. 
Regional-level daily hospitalization census data was 
obtained from Covid Patient Hospitalization Surveillance,  
posted publicly at 
\url{https://github.com/agb85/covid-19/tree/master/Regional%20Models}.
Parameter Values were obtained from previous modeling works 
\cite{AB-EC-DG-etal-b:21,AB-EC-DG-etal:21} and are reported in Table \ref{tab:modelParameters} and Supplementary Information.

Software that was custom-developed as part of our methods is 
available at the repository: 
\url{https://github.com/gianlucaBi/safe_levels_NPIs}.

\bibliography{alias,sample,GB}

\section*{Author contributions statement}

E.D., A.B., and G.B. designed the research; G.B. and E.D. 
performed research; G.B., E.D., and A.B. wrote the paper; J.I.P 
contributed with analytic development; E.J.C. supervised the 
research; D.J. provided data.  All authors reviewed the manuscript. 

\section*{Additional information}

\noindent
\textbf{Competing interests.} The authors declare no affiliation or involvement in any organization or entity with a financial or non-financial interest in the subject matter or materials discussed in this manuscript.

\begin{table}[h!]
\centering
 \begin{tabular}{||c c c c||} 
 \hline
 Symbol & Value & Description & Source\\ [0.5ex] 
 \hline\hline
 $\beta$ & 0.58 & transmission rate & fitted \\  \hline
$\theta$ & 0.77 & probability of vaccinating an individual in compartment 
$s$ & \cite{AB-EC-DG-etal-b:21}\\ 
\hline
$\delta$ & 0.02965/365 & daily death/birth rate 
& \cite{AB-EC-DG-etal-b:21}\\
\hline
$\sigma$ & 1/365 & 1/duration of natural immunity & \cite{AB-EC-DG-etal-b:21}\\
\hline
$\eta$ & 1/730 & 1/duration of vaccine immunity & \cite{AB-EC-DG-etal-b:21}\\
\hline
$\epsilon$ & 1/4.2 & 1/latency period & \cite{AB-EC-DG-etal-b:21}\\
\hline
$\gamma$ & 1/9 & rate of recovery & \cite{AB-EC-DG-etal-b:21}\\
\hline
$\kappa^{i \rightarrow h}$ & 0.0143762 & probability of hospitalization after infection & fitted \\
\hline
$\kappa^{i \rightarrow d}$ & 0.00262289 & probability of death after infection & \cite{AB-EC-DG-etal-b:21}\\
\hline
$\kappa^{h \rightarrow d}$ & 0.099204 & probability of death after hospitalization & \cite{AB-EC-DG-etal-b:21}\\
\hline
$\rho$ & $1/7.489$ & 1/hospitalization period & \cite{AB-EC-DG-etal-b:21}\\
\hline
$y$ & 15,000--25,000 & vaccination rate &\\
\hline
$\nu$ & 0.81 & vaccination efficacy & \cite{AB-EC-DG-etal-b:21} \\
\hline
$\sbs{N}{pop}$ & $5840795$ & state population size, CO, USA & \\
\hline
\hline
$s(0)$ & 1/1.47 & fraction of Susceptible on 03/01/21 & \cite{AB-EC-etal:21b}\\
\hline
$e(0)$ & 1/546 & fraction of Exposed on 03/01/21 & \cite{AB-EC-etal:21b}\\
\hline
$\imath(0)$ & 1/216 & fraction of Infectious on 03/01/21  & \cite{AB-EC-etal:21b}\\
\hline
$h(0)$ & 1/15936 & fraction of Hospitalized on 03/01/21 & \cite{AB-EC-etal:21b}\\
\hline
$r(0)$ & 1/4.2136 & fraction of Recovered on 03/01/21 & \cite{AB-EC-etal:21b}\\
\hline
$v(0)$ & 1/13.1 & fraction of Vaccinated on 03/01/21 & \cite{AB-EC-etal:21b}\\
\hline
$u(0)$ & 0.21 & level of lockdown on 03/01/21 & \cite{AB-EC-DG-etal-b:21}\\
\hline
\end{tabular}
\caption{Model parameters resulting from the model fitting phase.}
\label{tab:modelParameters}
\end{table}

\end{document}